\documentclass[12pt]{amsart}
\usepackage[]{amscd,amsfonts,amsmath,amsxtra,amssymb}
\usepackage{hyperref}
\usepackage{graphics}
\usepackage[all]{xy}
\newtheorem{sub}{}[section]
\newtheorem{subsub}{}[sub]

\def\EEnd{\mathop{\mathcal End}\nolimits}

\def\QQuot{\mathop{\mathcal Quot}\nolimits}

\def\hfl#1#2{\smash{\mathop{\ \hbox to 12mm{\rightarrowfill}}
\limits^{\scriptstyle#1}_{\scriptstyle#2} \ }}
\def\hflb#1#2{\smash{\mathop{\hbox to 12mm{\leftarrowfill}}
\limits^{\scriptstyle#1}_{\scriptstyle#2}}}

\def\og{\leavevmode\raise.3ex\hbox{$\scriptscriptstyle\langle\!\langle$}}
\def\fg{\leavevmode\raise.3ex\hbox{$\scriptscriptstyle\,\rangle\!\rangle$}}

\def\Ssect#1#2{\pagebreak[3]\begin{sub}\label{#2}{\sc\small\small
#1}\rm\medskip}

\def\xmat#1{\[\xymatrix{#1}\]}
\def\flinc{\ar@{^{(}->}}
\def\fleq{\ar@{=}}
\def\flon{\ar@{->>}}
\def\fmaps{\ar@{|-{>}}}

\def\BB{\mathop{\mathcal B}\nolimits}

\def\DD{\mathop{\mathcal D}\nolimits}

\def\SS{\mathop{\mathcal S}\nolimits}

\def\WW{\mathop{\mathcal W}\nolimits}

\def\ZZ{\mathop{\mathcal Z}\nolimits}

\newcommand{\N}{{\mathbb N}}

\newcommand{\Z}{{\mathbb Z}}

\newcommand{\C}{{\mathbb C}}

\renewcommand{\P}{{\mathbb P}}

\newcommand{\kf}{{\mathcal F}}

\newcommand{\kj}{{\mathcal J}}

\newcommand{\ko}{{\mathcal O}}

\newcommand{\kq}{{\mathcal Q}}

\newcommand{\ku}{{\mathcal U}}
\newcommand{\kv}{{\mathcal V}}

\newcommand\bigzero{\makebox(0,0){\text{\huge0}}}

\begin{document}

\def\refname{R\'ef\'erences}
\def\contentsname{ }
\def\proofname{D\'emonstration}
\def\abstractname{R\'esum\'e}

\author{Mohamed Bahtiti}
\address{Institut de Math\'{e}matiques de Jussieu,
 Case 247, 4 place Jussieu,\newline
F-75252 Paris, France}
\email{mohamed.bahtiti@imj-prg.fr }
\title{ Fibr\'e de Tango pond\'er\'e g\'en\'eralis\'e \mbox{de rang $n-1$~sur~l'espace $\P^{n}$} }
 
\maketitle

\begin{abstract}
Nous \'etudions dans cet article une nouvelle famille de fibr\'es vectoriels alg\'ebriques \mbox{stables} de rang $n-1$ sur l'espace projectif complexe $\P^{n}$ dont les fibr\'es de Tango pond\'er\'es de Cascini \cite{ca} font \mbox{partie}. Nous montrons que cette famille est invariante par rapport aux d\'eformations \mbox{miniversales}. 
$$$$

{\scriptsize ABSTRACT}. We study in this paper a new family of stable algebraic vector bundles of rank $ n-1 $ on the complex projective space $\P^{n}$ whose weighted Tango bundles of Cascini \cite{ca} belongs to. We show that these bundles are invariant under a miniversal deformation.
\end{abstract}

{\scriptsize {\em Date}: August, 2015.}\\
{\scriptsize 2010 Mathematics Subject Classification. 14D20, 14J60, 14F05, 14D15.}\\
{\scriptsize Mots-cl\'es. fibr\'e de Tango, stabilit\'e, d\'eformation miniversale, espace de Kuranishi.}\\
{\scriptsize Key words. Tango bundles, stability, miniversal deformation, Kuranishi space}.

\vspace{1.8cm}

\section{\large\bf Introduction}
 
\vspace{1cm}

Les fibr\'es vectoriels alg\'ebriques non-d\'ecomposables connus de rang $n-1$ sur l'espace projectif complexe $\P^{n}$ pour $n \geq 6$ sont rares. Les familles de fibr\'es vectoriels connues sont seulement la famille de fibr\'es instantons de rang $n-1$ pour $n$ impair \cite{ok-sp} et celle de fibr\'es de Tango pond\'er\'es de rang $n-1$ \cite{ca}. 

La famille de fibr\'es de Tango pr\'esente un sujet int\'eressant dans la g\'eom\'etrie alg\'ebrique. Cette famille de fibr\'es a \'et\'e construite sur $\P^{n}$ par Tango \cite{ta1}. Horrocks \cite{ho} a introduit une technique de construction de nouveau fibr\'e \`a partir d'un ancien fibr\'e muni d'une action de $\C^{*}$. Cette technique a \'et\'e appel\'ee l'image inverse g\'en\'eralis\'ee qui a \'et\'e \'etudi\'ee attentivement par Ancona et Ottaviani \cite{an-ot93}. En utilisant cette technique Cascini \cite{ca} a g\'en\'eralis\'e le fibr\'e de Tango, qui est $SL(2)$-invariant, \`a un fibr\'e de Tango pond\'er\'e. 

Dans cet article nous nous int\'eressons en particulier \`a la g\'en\'eralisation du fibr\'e de Tango qui est $\C^{*}$-invariant. Plus pr\'ecis\'ement, soient $i,n, \alpha ,\gamma \in \N$ et $\beta \in \Z$ tels que $n>2$, $\gamma >0$, $\alpha \geq \beta$, $\alpha +\beta\geq 0$ et $\gamma + n\alpha+i(\beta -\alpha)> 0$ pour $0\leq i\leq n$. Soient $Q$ le fibr\'e de quotient et $F(W)$ le fibr\'e de Tango sur $\P(\SS^{n}U)$ pour un $  W  \in \WW$, et $D$ comme dans le th\'eor\`eme \ref{3.4}. Alors le fibr\'e $Q$ (resp. $F(W)$) a une image invers\'ee g\'en\'eralis\'ee $Q_{\gamma,\alpha,\beta}$ (resp. $F_{\gamma,\alpha,\beta}$) d\'efinie par la suite exacte suivante
$$ 0\longrightarrow \ko_{\P^{n}}(-\gamma)\stackrel{}{\longrightarrow} \bigoplus_{i=0}^{n}\ko_{\P^{n}}(n\alpha+i(\beta -\alpha)) \longrightarrow \kq_{\gamma,\alpha,\beta} \longrightarrow 0$$

$$\text{ (resp. } 0\longrightarrow \kq_{\gamma,\alpha,\beta}(-\gamma) \stackrel{}{\longrightarrow}\bigoplus_{k=1}^{2n-1}\ko_{\P^{n}}(2n\alpha+k(\beta -\alpha)) \longrightarrow  \kf_{\gamma,\alpha,\beta}(\gamma) \longrightarrow 0),$$

o\`u $\kf_{\gamma,\alpha,\beta}(\gamma):=F_{\gamma,\alpha,\beta}(-2\gamma)$, $\kq_{\gamma,\alpha,\beta}:=Q_{\gamma,\alpha,\beta}(-\gamma)$.  On appelle le fibr\'e $\kq_{\gamma,\alpha,\beta}$ {\em le fibr\'e de quotient 
pond\'er\'e par les poids $\gamma,\alpha,\beta$, provenant d'une image inverse g\'en\'eralis\'ee sur $\P^{n}$}. On appelle le fibr\'e $\kf_{\gamma,\alpha,\beta} $ {\em le fibr\'e de Tango pond\'er\'e par les poids $\gamma,\alpha,\beta$, provenant d'une image inverse g\'en\'eralis\'ee sur $\P^{n}$}. En particulier le fibr\'e $\kf_{\gamma,\alpha,-\alpha}$ pond\'er\'e  par les poids $\gamma,\alpha ,-\alpha$ est le fibr\'e pond\'er\'e par les poids $\gamma,\alpha$ de Cascini \cite{ca}.

Le fibr\'e $\kf_{\gamma,\alpha,\beta}$ sur $\P^{n}$ v\'erifie les conditions suivantes (th\'eor\`eme \ref{3.9})

1- Si on a $\gamma > 2n\alpha+(\beta-\alpha)$, alors $\kf_{\gamma,\alpha,\beta}$ est stable.

2- Soit $\gamma> n\alpha$. Si $\kf_{\gamma,\alpha,\beta}$ est stable, alors on a $\gamma > 2n\alpha+(\beta-\alpha)$.

Les d\'eformations miniversales d'un tel fibr\'e $\kf_{\gamma,\alpha,\beta}$ sont encore des fibr\'es de Tango pond\'er\'es par les poids $\gamma,\alpha ,\beta$ sur $\P^{n}$ et l'espace de Kuranishi du fibr\'e $\kf_{\gamma,\alpha,\beta}$ est lisse au point correspondant de $\kf_{\gamma,\alpha,\beta}$ (th\'eor\`eme \ref{4.8}). 

Je tiens \`a exprimer ma gratitude au directeur de ma th\`ese M. J-M. Dr\'ezet et au professeur \mbox{G. Ottaviani} pour nos discussions utiles. Je remercie \'egalement toutes les personnes qui ont contribu\'e \`a m'aider \`a r\'ealiser mes travaux. Cet article fait partie de ma th\`ese.

\newpage

\section{\Large\bf Pr\'eliminaires} 
 
\vspace{1cm}

\subsection{\bf Remarque}\label{2.1} Si $D$ est un espace vectoriel de dimension $1$ et $q$ un entier, on note
$$D^{q}=\left\{
\begin{array}{ccc}

D\otimes D\otimes \ldots  \otimes D\otimes D \hspace{0.2 cm}\text{(q fois)}& :& q>0\\

\C & :& q=0\\

D^{*} \otimes D^{*} \otimes \ldots  \otimes D^{*} \otimes D^{*}\hspace{0.2 cm} \text{(-q fois)}& :& q<0.\\
 \end{array}
\right.
$$

et 
$$\SS^{q} D=\left\{
\begin{array}{ccc}

\SS^{q}D & :& q>0\\

\C & :& q=0\\

\SS^{-q} D^{*} & :& q<0.\\
 \end{array}
\right.
$$

M\^eme chose pour un fibr\'e vectoriel sur une vari\'et\'e $X$. On a donc, pour tout entier $q$ et $x=\C.v\in P^{n}= P(V)$,
$$(\ko_{\P(V)}(q))_{x} = x^{-q}=(\C.v)^{-q} .$$

On peut dire que $v^{*}:=v^{-1}$ est le vecteur dual de $ v$. On peut donc d\'efinir $ v^{-q}\in(\C.v)^{-q}\subset \SS^{-q}V$ pour tout entier $q$.
 
\subsection{D\'efinition de l'image invers\'ee g\'en\'eralis\'ee d'un fibr\'e (transformation de \mbox{Horrocks} \cite{ho})}\label{2.2}

Soient $V$ un espace vectoriel complexe de dimension $n+1$ et $\P^{n}=\P(V)$ l'espace projectif complexe associ\'e dont les points sont les droites de $V$. Soient
$$\eta:V\setminus \lbrace 0\rbrace \longrightarrow \P(V) $$

la projection, et $T$ une $\C^{*}$-action triviale (la multiplication usuelle) sur $V\setminus \lbrace 0\rbrace$ 
$$T:\C^{*}\times  V\setminus \lbrace 0\rbrace   \longrightarrow V\setminus \lbrace 0\rbrace$$
 $$(t,v)\longmapsto t.v=tv.$$

L'action $T$ induit une action triviale de $\C^{*}$ sur $\P(V)$ telle que $\eta$ est $\C^{*}$-\'equivariant et que $$\P(V)\simeq \left( V\setminus \lbrace 0\rbrace\right) \diagup \C^{*}     .$$

Soit $\kf\kv(\P(V))$ la cat\'egorie de fibr\'es vectoriels sur $\P(V)$. Soit $\kf\kv(V\setminus \lbrace 0\rbrace,T)$ la cat\'egorie de fibr\'es vectoriels sur $V\setminus \lbrace 0\rbrace$ qui sont $\C^{*}$-invariants au-dessus de l'action de $T$ sur $V\setminus \lbrace 0\rbrace$, ses morphismes \'etant des morphismes $\C^{*}$-\'equivariants des fibr\'es (les morphismes \'etant compatibles avec l'action $T$). Pour tout fibr\'e $E\in \kf\kv(\P(V))$, pour tout $t\in \C^{*} $ et $v\in V\setminus \lbrace 0\rbrace$, on a
$$(\eta ^{*}E)_{v}= E_{\eta (v)}\simeq E_{t.\eta (v)}=E_{\eta (T(t).v)}=(\eta ^{*}E)_{T(t).v}. $$ 

Donc on obtient un foncteur de cat\'egories
$$\kf\kv(\P(V))\stackrel{ \eta ^{*}(\bullet) }{ \longrightarrow }\kf\kv(V\setminus \lbrace 0\rbrace,T)$$
$$ E\longmapsto \eta^{*} E.$$

Par cons\'equent, pour tout fibr\'e $F\in  \kf\kv(V\setminus \lbrace 0\rbrace,T) $, il existe un fibr\'e $E \in \kf\kv(\P(V)) $ tel que l'on ait un $\C^{*}$-isomorphisme $F\simeq \eta^{*} E$. Par exemple, pour tout entier $q$ et pour tout $ v\in V \setminus \lbrace 0\rbrace $, on a un isomorphisme canonique
$$\C \stackrel{\simeq}{\longrightarrow} \eta^{*}(\ko_{\P(V)} (q))_{v} $$
$$t\longmapsto t.v^{-q}.$$

Cela d\'efinit un isomorphisme
$$\ko_{V\setminus \lbrace 0\rbrace} \stackrel{\simeq}{\longrightarrow} \eta^{*}(\ko_{\P(V)} (q)).$$

Le fibr\'e $\eta^{*}(\ko_{\P(V)} (q))$ a une action canonique de $\C^{*}$ compatible avec l'action de ce groupe sur $V\setminus \lbrace 0\rbrace$. Compte tenu de l'isomorphisme pr\'ec\'edent, c'est une action de $\C^{*}$ sur $\ko_{V\setminus \lbrace 0\rbrace}$.
Cette action est la suivante
$$\C^{*}\times \ko_{V\setminus \lbrace 0\rbrace} =\C^{*}\times (  V\setminus \lbrace 0\rbrace \times \C) \longrightarrow V\setminus \lbrace 0\rbrace \times \C $$
$$\hspace{4 cm}(t,(u,a))\longmapsto (tu,t^{q}a)$$

Autrement dit, cette action est la multiplication de l'action triviale de $\C^{*}$ sur $\ko_{V\setminus \lbrace 0\rbrace}$ par le~caract\`ere $t^{q}$. Donc $\eta^{*}(\ko_{\P(V)} (q))$ est le fibr\'e trivial sur $V\setminus \lbrace 0\rbrace$ muni de l'action pr\'ec\'edente.

Cette correspondance est compatible avec les op\'erations habituelles sur les fibr\'es. Par exemple, si on consid\`ere  que $F$ (resp. $F^{*}$) est $\C^{*}$-invariant au-dessus de l'action $T$ (resp. $\C^{*}$-invariant au-dessus de l'action $\widehat{T}$ qui est l'action duale de $T$), alors on a $\eta^{*}(E^{*})=F^{*}$. M\^eme chose pour les produits tensoriels (resp. ext\'erieurs, sym\'etriques), la somme directe, $Hom(,)$ et une suite exacte (resp. une monade) de fibr\'es vectoriels de trois termes.

Soient $g_{0},g_{1},\ldots ,g_{n}$ des polyn\^omes homog\`enes de degr\'es $d_{0}\geq d_{1} \geq  \ldots \geq d_{n} $ respectivement sans z\'ero commun sur $\P(V_{2})$. On a l'application surjective
$$\omega:= (g_{0},g_{1},\ldots ,g_{n}): V_{1}\setminus  \lbrace 0 \rbrace \longrightarrow V_{2}\setminus  \lbrace 0 \rbrace \hspace{2 cm}$$
$$\hspace{4 cm}v\longmapsto (g_{0}(v),g_{1}(v), \ldots ,g_{n}(v))$$

o\`u $V_{1}=\C^{n+1}$ et $V_{2}$ est un $\C$-espace vectoriel de dimension $n+1$. Soit \mbox{$\eta_{i}: V_{i}\setminus  \lbrace 0 \rbrace \longrightarrow \P( V_{i})$} la projection pour $i=1,2$. On consid\`ere l'action de $\C^{*}$ sur $V_{2}$
$$\sigma :\C^{*} \longrightarrow GL(V_{2}) $$
$$ t\longmapsto \sigma(t)$$
 
o\`u 
$$\sigma(t)=\left(
\begin{array}{ccccccccc}
 t^{d_{0}} &  &   & & &            &   \\
     & t^{d_{1}} &  &   & &   &\bigzero \\
  &  & \ddots &   &    &   &   \\
\bigzero &  &   &    &  t^{d_{n-1}}&  &   \\
 &  &   &    &    &  &      t^{d_{n}}\\
\end{array}
\right),$$

et on consid\`ere une $\C^{*}$-action $T$ qui est la multiplication usuelle sur $V_{1}\setminus \lbrace 0\rbrace$ 
$$T:\C^{*}\times  V_{1}\setminus \lbrace 0\rbrace   \longrightarrow V_{1}\setminus \lbrace 0\rbrace$$
$$(t,v)\longmapsto T(t,v)=t.v $$

de telle sorte que $\eta_{i}$ est un $\C^{*}$-morphisme. Alors $\omega $ est une $\C^{*}$-application par rapport \`a ces deux actions. L'action $\sigma$ induit une action $\overline{ \sigma}\in PGL(V_{2}) $ de $\C^{*}$ sur $\P(V_{2})$ et l'action $T$ induit une action triviale de $\C^{*}$ sur $\P(V_{1})$. On obtient donc le $\C^{*}$-diagramme suivant
\xmat{ V_{1} \setminus  \lbrace 0 \rbrace \ar[r]^{\omega } \ar[d]_{\eta_{1}} 
&  V_{2} \setminus \lbrace 0 \rbrace \ar[d]^{\eta_{2}} \\
\P^{n}=\P(V_{1})  &  \P^{n}=\P(V_{2}) }

Soit $F$ un fibr\'e vectoriel sur $\P(V_{2})$ qui est $\C^{*}$-invariant au-dessus de l'action $\overline{ \sigma}$. Donc $ \eta_{2}^{*}F$ est $\C^{*}$-invariant au-dessus de l'action $\sigma$. Alors $\omega^{*}\eta_{2}^{*}F$ est $\C^{*}$-invariant au-dessus de l'action usuelle de $\C^{*}$ sur $V_{1}\setminus \lbrace 0\rbrace$. Autrement dit, pour tout $v\in V_{1} \setminus  \lbrace 0 \rbrace $ et $t\in \C^{*}$, on a
$$ (\omega ^{*}\eta_{2} ^{*}F)_{v}= F_{\eta_{2} (\omega(v))}\simeq F_{\overline{\sigma(t)}.\eta_{2} (\omega(v))}=F_{\eta_{2} (\sigma(t).\omega(v))}=(\eta_{2} ^{*}F)_{\sigma(t).\omega(v)}=(\eta_{2} ^{*}F)_{\omega(t.v)}=(\omega ^{*}\eta_{2} ^{*}F)_{t.v}. $$

Alors il existe un fibr\'e vectoriel $F_{f_{0},f_{1},\ldots ,f_{n}}$ sur $\P(V_{1})$ tel que l'on ait un $\C^{*}$-isomorphisme 
$$\omega^{*}\eta_{2}^{*}F\simeq \eta_{1}^{*}F_{f_{0},f_{1},\ldots ,f_{n}}.$$

On appelle $F_{1}:= F_{f_{0},f_{1},\ldots ,f_{n}}$ l'image inverse g\'en\'eralis\'ee de $F$. Donc on a le foncteur 
\xmat{  \kf\kv(\P(V_{2}),\overline{ \sigma})   \ar[rr]^{{\bf Iminvg}} &  &  \kf\kv(\P(V_{1}))  \\
 }
\xmat{ F \ar@{|->}[rr] & &  F_{1} }

et on a
\xmat{   \kf\kv(\P(V_{2}),\overline{ \sigma})   \ar[rr]^{\omega^{*}\eta_{2}^{*}(\bullet)}\ar[rd]^{{\bf Iminvg}} &  
&  \kf\kv(V_{1}\setminus \lbrace 0\rbrace,T)  \\
&    \kf\kv(\P(V_{1})) \ar[ru]^{\eta_{1}^{*}(\bullet)}}

 o\`u $\kf\kv(\P(V_{2}),\overline{ \sigma})$ est la cat\'egorie de fibr\'es vectoriels sur $\P(V_{2})$ qui sont $\C^{*}$-invariants au-dessus de l'action $\overline{ \sigma}$, et les morphismes sont des morphismes $\C^{*}$-\'equivariants des fibr\'es (les morphismes \'etant compatibles avec l'action $\overline{ \sigma}$). 

Cette {\em transformation de Horrocks} ${\bf Iminvg}$ est compatible avec les op\'erations habituelles sur les fibr\'es. Par exemple, si on consid\`ere que $F$ (resp. $F^{*}$) est $\C^{*}$-invariant au-dessus de l'action $\overline{ \sigma}$ \mbox{(resp. $\C^{*}$-invariant au-dessus de l'action }$\widehat{\overline{ \sigma}}$ qui est l'action duale de $\overline{ \sigma}$), alors on a ${\bf Iminvg}(F^{*})={\bf Iminvg}(F)^{*}$. M\^eme chose pour les produits tensoriels (resp. ext\'erieurs, sym\'etriques), la somme directe, $Hom(,)$ et une suite exacte (resp. une monade) de fibr\'es vectoriels de trois termes.

\subsection{\bf Proposition}\label{2.3}
{\em On consid\`ere les m\^emes notations de la d\'efinition \ref{2.2}. Soient $E $ un \mbox{$C^*$-fibr\'e} quelconque sur $\P(V_{2})$, $s : C^{*} \times E \longrightarrow E$ son $C^*$-action au-dessus de l'action $\overline{ \sigma}$ et $q$ un entier. On en d\'eduit une nouvelle action de $C^*$ sur $E$, pour tout $x\in \P(V_{2})$,
$$s_{q}{}_{,x} : C^{*}\times E_{x} \longrightarrow E_{x}$$
        $$\hspace{2.9cm}(t,u) \longmapsto t^{q}.s_{x}(t,u)$$

(multiplication de l'action par le caract\`ere $t^{q}$). On note le $C^*$-fibr\'e obtenu $E^{(q)}=E\otimes\ko_{\P(V_{2})}^{(q)} $. Alors on a

1- ${\bf Iminvg}(\ko_{\P(V_{2})}^{(q)}(k))=\ko_{\P(V_{1})}(q)$, pour tout entier $k$.

2- ${\bf Iminvg}(E^{(q)}) = ({\bf Iminvg}(E))(q)$.} 

\begin{proof}
1- On a ${\bf Iminvg} (\ko_{\P(V_{2})}^{(q)}(k))= \ko_{\P(V_{1})}(d)$, o\`u $d$ est un entier. C'est-\`a-dire 
$$ \omega^{*}\eta_{2}^{*} (\ko_{\P(V_{2})}^{(q)}(k))\simeq \eta_{1}^{*}\ko_{\P(V_{1})}(d) .$$

On a, pour tout $x=\C.v_{2}\in \P(V_{2})$, un isomorphisme canonique
$$ \mu_{v_{2}}:\ko_{V_{2}\setminus \lbrace 0\rbrace,v_{2}}\simeq \C \stackrel{. v_{2}^{-k} }{\longrightarrow} \eta_{2}^{*}(\ko_{\P(V_{2})}^{(q)}(k))_{v_{2}}\simeq (\C v_{2})^{-k} \subset \SS^{-k}V_{2} $$
\xmat{ a \ar@{|->}[rrrr] & &&&   a.v_{2}^{-k} .}

On a aussi, pour tout $x=\C.v_{1}\in \P(V_{1})$, un isomorphisme canonique
$$ \overline{\mu}_{v_{1}}:\ko_{V_{1}\setminus \lbrace 0\rbrace,v_{1}}\simeq \C \stackrel{.v_{1}^{-d} }{\longrightarrow} \eta_{1}^{*}(\ko_{\P(V_{1})}(d))_{v_{1}}\simeq (\C.v_{1})^{-d}\subset \SS^{-d}V_{1} $$
\xmat{ a \ar@{|->}[rrrr] & &&&   a.v_{1}^{-d} .}

On a, pour tout $v\in V_{1}$, un isomorphisme  
$$(\C \omega(v))^{-k} \simeq (\omega^{*}\eta_{2}^{*} \ko_{\P(V_{2})}^{(q)}(k))_{v}\simeq (\eta_{1}^{*}\ko_{\P(V_{1})}(d))_{v}\simeq (\C.v)^{-d}$$
\xmat{ a.(\omega(v))^{-k} \ar@{|->}[rrrrrr] & &&&& &  a.(v)^{-d} .}

Alors on obtient le diagramme commutatif suivant

 \xmat{& & \ar@/_4pc/[4,0]_{v^{-d}} \ar[dd]^{(\omega (v))^{-k}} \C^{*} \ar[rr]^{t^{q}}& &  \ar[dd]_{(\omega (t.v))^{-k}} \C^{*} \ar@/^4pc/[4,0]^{(t.v)^{-d}} & &\\ \\
& &\ar@{=}[dd]_{}(\C \omega(v))^{-k} \ar[rr]^{s_{q,}{}_{v}} & &   \ar@{=}[dd]^{} (\C \omega(t.v))^{-k}  & & \\ \\
 & & (\C.v)^{-d} \ar[rr] & &  (\C.v)^{-d} \simeq (\C.(t.v))^{-d} & &  } 

tel que $a.(v)^{-d} = a. t^{q}(t.v)^{-d}:=a. t^{q-d}(v)^{-d}$, o\`u $a\in \C^{*}$. Par cons\'equent $d=q$ et on obtient 
$${\bf Iminvg}(\ko_{\P(V_{2})}^{(q)}(k))=\ko_{\P(V_{1})}(q).$$

2- Comme $E^{(q)}=E\otimes\ko_{\P(V_{2})}^{(q)} $ et que le foncteur ${\bf Iminvg}$ respecte le produit \mbox{tensoriel}, on obtient
$${\bf Iminvg}(E^{(q)})={\bf Iminvg}(E)\otimes {\bf Iminvg}(\ko_{\P(V_{2})}^{(q)}).$$

D'apr\`es (1), on a \mbox{${\bf Iminvg}(\ko_{\P(V_{2})}^{(q)})=\ko_{\P(V_{1})}(q)$. On en d\'eduit donc} $${\bf Iminvg}(E^{(q)})={\bf Iminvg}(E)(q).$$

\end{proof}

\newpage

\section{\Large\bf  Fibr\'e de Tango pond\'er\'e g\'en\'eralis\'e } 
 
\vspace{1cm}

Le fibr\'e de Tango, qui est $SL_{2}(\C)$-invariant, et son fibr\'e de Tango pond\'er\'e provenant d'une image inverse g\'en\'eralis\'ee sur $\P^{n}$ ont d\'ej\`a \'et\'e trait\'es dans l'article \cite{ca}. Nous allons traiter le fibr\'e de Tango, qui est $\C^{*}$-invariant, et son fibr\'e de Tango pond\'er\'e provenant d'une image inverse g\'en\'eralis\'ee sur $\P^{n}$.

\subsection{\bf D\'efinition} (Jaczewski, Szurek, Wisniewski \cite{ja-sz-wi} et Tango \cite{ta1}). \label{3.1} 
Soient $V$ un $\C$-espace vectoriel de dimension $dim(V)=n+1$, et $\P^{n}=P(V)$ 
l'espace projectif complexe associ\'e \`a l'espace $V$. Soit $W \subset \bigwedge^{2}V=(H^{0}(Q^{*}(1)))^{*} $ un sous-espace vectoriel tel que
$$(*)\left\{
\begin{array}{ccc}

\hspace{0.2cm} - \ dim(W)= \left(\begin{array}{c}n+1\\2\end{array}\right)- 2n+1  \hspace{5.8cm}\\

- \ W \hspace{0.2cm} \text{ne contient pas d'\'el\'ement d\'ecomposable non nul de}\hspace{0.2cm} \bigwedge^{2}V.

 \end{array}
\right.
$$

En utilisant la suite exacte suivante
$$ 0\longrightarrow \ko_{\P^{n}}(-1)\stackrel{g}{\longrightarrow} V\otimes \ko_{\P^{n}}\stackrel{\mu}{\longrightarrow} Q\longrightarrow 0,$$

on obtient la r\'esolution suivante
$$0 \longrightarrow \ko_{\P^{n}}(-2) \stackrel{g\otimes I_{\ko_{\P^{n}}(-1)}}{\longrightarrow}    V \otimes \ko_{\P^{n}}(-1)  \stackrel{I_{V } \bigwedge g}{\longrightarrow} \bigwedge^{2}V\otimes \ko_{\P^{n}} \stackrel{\bigwedge^{2}\mu }{\longrightarrow} \bigwedge^{2} Q \longrightarrow  0. $$ 

Donc on en d\'eduit la suite exacte suivante
$$ 0\longrightarrow Q(-1)\stackrel{\beta}{\longrightarrow} \bigwedge^{2} V\otimes \ko_{\P^{n}}\stackrel{\bigwedge^{2}\mu }{\longrightarrow} \bigwedge^{2} Q\longrightarrow 0,$$

o\`u $I_{V} \wedge g =\beta \circ (\mu\otimes I_{\ko_{\P^{n}}(-1)} ) $. On a le morphisme d'\'evaluation du fibr\'e $Q^{*}(1)$ 
$$ev_{Q^{*}(1)}:\ko_{\P^{n}}\otimes  \bigwedge^{2} V \longrightarrow Q^{*}(1). $$

Il en d\'ecoule que $\beta =^{T}ev_{Q^{*}(1)}$. Pour tout $x=\C.v_{0}\in \P^{n}$ et $v_{0}\in V$, on a 

\xmat{ 0 \ar[rr] && W  \ar[rr]  &&\bigwedge^{2}V \ar[rr]^{q} &&(\bigwedge^{2}V)\diagup W \ar[rr]&&0 \\
  &&  && Q_{x}(-1)\ar[rru]_{(\varpi_{W})_{x}}\ar[u]^{\beta_{x}}\\
  && &&\ar[u] 0 }

L'application $(\varpi_{W})_{x}= q\circ \beta_{x}$ est injective car: 

soit $(\varpi_{W})_{x}(a)=0$, pour tout $a \in Q_{x}(-1)$, on obtient que $\beta_{x}(a)\in W $. Il existe \'egalement $v\in V$ tel que 
$$a=(\mu\otimes I_{\ko_{\P^{n}}(-1)} )_{x} (v\otimes x)=\mu_{x}(v)\otimes v_{0}.$$

Donc on a
$$ (\beta \circ (\mu\otimes I_{\ko_{\P^{n}}(-1)} ))_{x}  (v\otimes x)=(I_{V}  \wedge g)_{x} (v\otimes x), $$
$$\beta_{x} (a)=v \wedge g_{x}(x)  =v \wedge  v_{0}.$$

Comme $W$ ne contient pas d'\'el\'ement d\'ecomposable non nul de $\bigwedge^{2}V$, alors on obtient $\beta_{x}(a)=0$ qui donne $a=0$.
 On d\'efinit {\em le fibr\'e de Tango $F(W)$ de rang $n-1$ sur $\P^{n}$} par la suite exacte suivante
$$0\longrightarrow Q(-1)\stackrel{\varpi_{W}}{\longrightarrow} \left(    (\bigwedge^{2}V)\diagup W \right)\otimes\ko_{\P^{n}} \longrightarrow  F(W)(1)   \longrightarrow 0.$$  

Sa premi\`ere classe de Chern est $c_{1}=2n$ et $ H^{0}( F(W)(1))=(\bigwedge^{2}V)\diagup W$. On a un carr\'e commutatif

\xmat{  \left((\bigwedge^{2}V)\diagup W \right)^{*}\otimes\ko_{\P^{n}} \ar[rr]^{^{T}\varpi_{W}} \ar@{^{(}->}[d]&& Q^{*}(1) \ar@{=}[d]\\
  (\bigwedge^{2}V)^{*}\otimes\ko_{\P^{n}}  \ar[rr]^{ev_{Q^{*}(1)}}  & &  Q^{*}(1).\\
}
  
On en d\'eduit imm\'ediatement que l'inclusion $^{T}q:\left((\bigwedge^{2}V)\diagup W \right)^{*}\subset (\bigwedge^{2}V)^{*} $ s'identifie \`a $H^{0}(^{T}\varpi_{W})$. De la troisi\`eme suite exacte, il d\'ecoule la suite de cohomologies suivante
$$H^{0}((F(W)(1))^{*})=Hom(F(W)(1),\ko_{\P^{n}}) \longrightarrow Hom(\left((\bigwedge^{2}V)\diagup W \right)\otimes\ko_{\P^{n}},\ko_{\P^{n}})\stackrel{H^{0}(^{T}\varpi_{W})}{\longrightarrow} $$
$$Hom(Q(-1),\ko_{\P^{n}}) \longrightarrow Ext^{1}(F(W)(1),\ko_{\P^{n}}) \longrightarrow 0.$$

Donc on a $H^{0}((F(W)(1))^{*})=0$ et $H^{1}((F(W)(1))^{*})=Ext^{1}(F(W)(1),\ko_{\P^{n}})=W^{*}$.

\subsection{\bf Proposition }\label{3.2}
{\em Soit $\bigwedge^{2}V=(H^{0}(Q^{*}(1)))^{*}$. On a les assertions suivantes

1- Soient $W_{1}, W_{2}$ des sous-espaces vectoriels de $\bigwedge^{2}V$ v\'erifiant la condition (*). Alors on obtient 
$F(W_{1})\simeq F(W_{2})$ si et seulement si on a $W_{1}= W_{2}$.

2- Soient $\sigma \in GL(V)$ et $\overline{\sigma } \in PGL(V) $ son \'el\'ement correspondant. Soit $W$ un sous-espace vectoriel de $\bigwedge^{2}V$ v\'erifiant la condition (*). Alors $\sigma^{-1}(W)$ est un sous-espace vectoriel de $\bigwedge^{2}V$ v\'erifiant la condition (*) et il existe un isomorphisme canonique
$$\Psi^{W}_{\sigma}:\overline{\sigma }^{*}(F(W)) \stackrel{\simeq}{\longrightarrow}F(\sigma^{-1}(W)).$$

Soient $\rho \in GL(V)$ et $\overline{\rho } \in PGL(V) $ son \'el\'ement correspondant, alors on a
$$\Psi^{W}_{\sigma \rho}=\Psi^{\sigma^{-1}(W)}_{\rho}\circ \overline{\rho}^{*}\Psi^{W}_{\sigma}.$$ 

En particulier, pour tout $t\in \C^{*}$, soient
 $$ \sigma(t)=\rho(t):=\left(
\begin{array}{ccccccccc}
 t^{a_{0}} &  &   & & &            &   \\
     & t^{a_{1}} &  &   & &   &\bigzero \\
  &  & \ddots &   &    &   &   \\
\bigzero &  &   &    &  t^{a_{n-1}}&  &   \\
 &  &   &    &    &  &      t^{a_{n}}\\
\end{array}
\right) \in GL(V)
$$

et $\overline{ \sigma(t)}, \overline{\rho(t)}\in PGL(V)$ leurs \'el\'ements correspondants tels que $ \sigma(t)(W)\subseteq W$. Alors il existe un isomorphisme canonique
$$s_{t}:\overline{ \sigma(t)}^{*}F(W)\simeq F(W)$$
  
tel que, pour $t_{1},t_{2}\in \C^{*}$, on a $ s_{t_{1}.t_{2}}=s_{t_{2}}\circ \overline{ \sigma(t_{2})}^{*}(s_{t_{1}})$. Autrement dit, on a le diagramme commutatif suivant

\xmat{   \overline{\sigma(t_{2})}^{*} \overline{\sigma(t_{1})}^{*}F(W)=\overline{\sigma( t_{1}.t_{2} )}^{*} F(W) \ar[rr]^{s_{t_{1}. t_{2}}}\ar[rd]_{\overline{\sigma(t_{2})}^{*}(s_{t_{1}) }} &  
&   F(W)  \\
&    \overline{\sigma(t_{2})}^{*} F(W)  \ar[ru]_{s_{t_{2}}}.}}

\begin{proof}

1- Pour $W_{1}, W_{2}$ des sous-espaces vectoriels de $\bigwedge^{2}V$ v\'erifiant la condition (*), on a les suites exactes suivantes
$$0\longrightarrow Q(-1)\stackrel{\varpi_{W_{1}} }{\longrightarrow}  \left(    (\bigwedge^{2}V)\diagup W_{1} \right)\otimes\ko_{\P^{n}} \stackrel{b_{1} }{\longrightarrow} F(W_{1})(1) \longrightarrow 0 $$ 

et 
$$0\longrightarrow Q(-1)\stackrel{\varpi_{W_{2}} }{\longrightarrow}  \left(    (\bigwedge^{2}V)\diagup W_{2} \right)\otimes\ko_{\P^{n}} \stackrel{b_{2} }{\longrightarrow} F(W_{2})(1) \longrightarrow 0. $$

Soit $\varphi :F(W_{1})(1)\simeq F(W_{2})(1)$. On d\'eduit de la deuxi\`eme suite
$$0\longrightarrow Hom( \left(    (\bigwedge^{2}V)\diagup W_{1} \right)\otimes\ko_{\P^{n}},Q(-1)  )\stackrel{}{\longrightarrow}$$
$$  Hom( \left(    (\bigwedge^{2}V)\diagup W_{1} \right)\otimes\ko_{\P^{n}}, \left(    (\bigwedge^{2}V)\diagup W_{2} \right)\otimes\ko_{\P^{n}}  )  \stackrel{ b_{2}\circ \bullet }{\longrightarrow}$$
$$ Hom( \left(    (\bigwedge^{2}V)\diagup W_{1} \right)\otimes\ko_{\P^{n}},  F(W_{2})(1)) \longrightarrow Ext^{1}(\left(    (\bigwedge^{2}V)\diagup W_{1} \right)\otimes\ko_{\P^{n}},Q(-1) )\longrightarrow 0. $$

Comme $H^{1}(Q(-1))=0$, alors le morphisme $ b_{2}\circ \bullet $ est surjectif. Donc pour le morphisme $\varphi\circ b_{1}$ il existe un morphisme 
$$H^{0}(\varphi)\otimes I_{\ko_{\P^{n}}}:  \left(    (\bigwedge^{2}V)\diagup W_{1} \right)\otimes\ko_{\P^{n}}  \longrightarrow  \left(    (\bigwedge^{2}V)\diagup W_{2} \right)\otimes\ko_{\P^{n}}$$

tel que $b_{2}\circ (H^{0}(\varphi)\otimes I_{\ko_{\P^{n}}}) =\varphi\circ b_{1}$. Comme $\varphi$ est un isomorphisme, alors $H^{0}(\varphi)\otimes I_{\ko_{\P^{n}}}$ est aussi un isomorphisme lequel d\'efinit un isomorphisme $\varphi_{2}: Q(-1)\longrightarrow Q(-1)$ qui est une homoth\'etie car le fibr\'e $Q$ est simple. Donc on obtient le diagramme commutatif suivant

\xmat{ 0 \ar[r] & Q(-1) \ar[rr]^{\varpi_{W_{1}}}\ar[d]_{\varphi_{2}} & & \left((\bigwedge^{2}V)\diagup W_{1} \right)\otimes\ko_{\P^{n}} \ar[rr] \ar[d]_{\wr}^{H^{0}(\varphi)\otimes I_{\ko_{\P^{n}}}}  && F(W_{1})(1)) \ar[r]\ar[d]_{\wr}^{\varphi}& 0\\
0 \ar[r] & Q(-1) \ar[rr]^{\varpi_{W_{2}}} &  &   \left( (\bigwedge^{2}V)\diagup W_{2} \right)\otimes\ko_{\P^{n}}  \ar[rr]  & &  F(W_{2})(1) \ar[r]& 0.\\
}

En consid\'erant la cohomologie de ce diagramme, on obtient 

\xmat{  Hom(\left((\bigwedge^{2}V)\diagup W_{2} \right)\otimes\ko_{\P^{n}},\ko_{\P^{n}})
\ar@{^{(}->}[rr]^{H^{0}(^{T}\varpi_{W_{2}})}\ar[d]_{\wr}^{^{T}H^{0}(\varphi)} & & Hom(Q(-1),\ko_{\P^{n}}) \ar@{->>}[r] \ar@{=}[d]  & Ext^{1}(F(W_{2})(1),\ko_{\P^{n}})\ar[d]_{\wr}^{^{T}H^{1}(\varphi)}\\
  Hom(\left((\bigwedge^{2}V)\diagup W_{1} \right)\otimes\ko_{\P^{n}},\ko_{\P^{n}})
\ar@{^{(}->}[rr]^{H^{0}(^{T}\varpi_{W_{1}})} & & Hom(Q(-1),\ko_{\P^{n}}) \ar@{->>}[r]   & Ext^{1}(F(W_{1})(1),\ko_{\P^{n}}) \\
}

Donc on obtient un diagramme commutatif

\xmat{ 0 \ar[r] & \left((\bigwedge^{2}V)\diagup W_{2} \right)^{*} \ar[rr]^{H^{0}(^{T}\varpi_{W_{2}})}\ar[d]^{^{T}H^{0}(\varphi)} & & (\bigwedge^{2}V)^{*} \ar[rr] \ar@{=}[d]  && (W_{2})^{*} \ar[r]\ar[d]^{^{T}H^{1}(\varphi)}& 0\\
0 \ar[r] & \left((\bigwedge^{2}V)\diagup W_{1} \right)^{*} \ar[rr]^{H^{0}(^{T}\varpi_{W_{1}})} &  &(\bigwedge^{2}V)^{*} \ar[rr]  & &  (W_{1})^{*} \ar[r]& 0.\\
}

Comme $H^{0}(^{T}\varpi_{W_{1}}), H^{0}(^{T}\varpi_{W_{2}})$ sont les inclusions naturelles, alors les
sous-espaces vectoriels $\left((\bigwedge^{2}V)\diagup W_{1} \right)^{*}$ et $\left((\bigwedge^{2}V)\diagup W_{2} \right)^{*}$ de $(\bigwedge^{2}V)^{*}$ sont \'egaux et donc $W1 = W2$.

2- Soient $\sigma \in GL(V)$ et $\overline{\sigma } \in PGL(V) $ son \'el\'ement correspondant. Alors on obtient que \mbox{$dim(W)=dim(\sigma^{-1}(W))$} et que $\sigma^{-1}(W)$ ne contient pas d'\'el\'ement d\'ecomposable non nul de $\bigwedge^{2}V $. Car si $\sigma^{-1}(W)$ contient $ \sigma^{-1}(y)\wedge \sigma^{-1}(z)=\sigma^{-1}(y\wedge z)$ qui est un \'el\'ement non nul, alors $W$ contient $y\wedge z$ qui est un \'el\'ement non nul; ce qui est une contradiction. On a la suite exacte suivante
$$0\longrightarrow \overline{\sigma }^{*}Q(-1)\stackrel{\overline{\sigma }^{*}\varpi_{W}}{\longrightarrow}  \overline{\sigma }^{*}\left( \left(    (\bigwedge^{2}V)\diagup W \right)\otimes\ko_{\P^{n}}\right) \stackrel{\overline{\sigma }^{*}b }{\longrightarrow} \overline{\sigma }^{*}F(W) \longrightarrow 0. $$

Comme on a, pour tout $x=\C.v \in \P^{n} $ o\`u $v\in V$,
$$\overline{\sigma }^{*}\left( \left((\bigwedge^{2}V)\diagup W \right)\otimes\ko_{\P^{n}}\right)_{x}=
\sigma. \left((\bigwedge^{2}V)\diagup W \right)\otimes\ko_{\P^{n}}{}_{,\overline{\sigma }(x)}$$
$$ =\sigma. \left((\bigwedge^{2}V)\diagup W \right) =(\bigwedge^{2}V) \diagup \sigma^{-1}(W) =\left(\left( (\bigwedge^{2}V) \diagup \sigma^{-1}(W)\right)\otimes\ko_{\P^{n}}\right)_{x},$$
  
alors on obtient 
$$0\longrightarrow \overline{\sigma }^{*}Q(-1)\stackrel{\overline{\sigma }^{*}\varpi_{W}}{\longrightarrow}  \left(    (\bigwedge^{2}V)\diagup \sigma^{-1}(W) \right)\otimes\ko_{\P^{n}} \stackrel{\overline{\sigma }^{*}b }{\longrightarrow} \overline{\sigma }^{*}F(W) \longrightarrow 0. $$ 

Comme le fibr\'e $Q$ est homog\`ene, alors on a  
$$0\longrightarrow Q(-1)\stackrel{\overline{\sigma }^{*}\varpi_{W} }{\longrightarrow}  \left(    (\bigwedge^{2}V)\diagup \sigma^{-1}(W) \right)\otimes\ko_{\P^{n}} \stackrel{\overline{\sigma }^{*}b }{\longrightarrow} \overline{\sigma }^{*}F(W) \longrightarrow 0. $$

Comme on a $\sigma^{*}\varpi_{W}=\varpi_{\sigma^{-1}(W)}$, alors on obtient

\xmat{ 0 \ar[r] & Q(-1) \ar[rrr]^{\overline{\sigma }^{*}\varpi_{W}}\ar[d]^{\wr}&& &
 \left((\bigwedge^{2}V)\diagup \sigma^{-1}(W) \right)\otimes\ko_{\P^{n}}\ar[r]^{ }\ar[d]^{\wr}  & \overline{\sigma }^{*}F(W) \ar[r]\ar[d]& 0\\
0 \ar[r] & Q(-1) \ar[rrr]^{\varpi_{\sigma^{-1}(W)}} & & &\left( (\bigwedge^{2}V)\diagup \sigma^{-1}(W) \right)\otimes\ko_{\P^{n}} \ar[r]^{ }  &   F(\sigma^{-1}(W)) \ar[r]& 0.\\
}

Ce qui implique qu'il existe un isomorphisme canonique
$$\Psi^{W}_{\sigma}:\overline{\sigma }^{*}(F(W)) \stackrel{\simeq}{\longrightarrow}F(\sigma^{-1}(W)).$$

Soient $\rho \in GL(V)$ et $\overline{\rho } \in PGL(V) $ son \'el\'ement correspondant. Alors on a un carr\'e commutatif

\xmat{  \overline{\rho }^{*}\overline{\sigma }^{*}F(W)=(\overline{\sigma }\circ\overline{\rho })^{*}F(W) \ar[rr]^{\Psi^{W}_{\sigma \circ\rho}} \ar[d]_{\overline{\rho }^{*}(\Psi^{W}_{\sigma})}&& F((\sigma \circ\rho )^{-1}(W)) \ar@{=}[d]\\
\overline{\rho }^{*} F(\sigma^{-1}(W))  \ar[rr]^{  \Psi^{\sigma^{-1}(W)}_{\rho}  }  & &  F(\rho^{-1}(\sigma^{-1}(W)))=F((\sigma \circ\rho )^{-1}(W)).}

Comme on a $ \sigma(t)(W)\subseteq W$ alors on en d\'eduit $ \sigma(t)^{-1}(W)= W$, pour tout $t\in \C^{*}$, et il existe un isomorphisme canonique
$$s_{t}:\overline{ \sigma(t)}^{*}F(W)\simeq F(W).$$
 
Comme $\sigma(t_{1}).\sigma(t_{2})=\sigma(t_{1}.t_{2})$, alors on a 
$$ s_{t_{2}.t_{1}}=s_{t_{2}}\circ \overline{ \sigma(t_{2})}^{*}(s_{t_{1}})$$

pour tout $t_{2},t_{1}\in \C^{*}$, tout en consid\'erant $\sigma (t) =\rho(t)$ dans le carr\'e commutatif pr\'ec\'edent. 

\end{proof} 

\subsection{\bf Remarque} \label{3.3} Soient $U$ un $\C$-espace vectoriel de dimension $2$, $\lbrace x,y\rbrace$ sa base et $n>2$ un entier. Soit 
$$\BB_{0}:=\lbrace v_{p}:=x^{n-p}y^{p}  ,  \hspace{0.2 cm} 0\leq p\leq n \rbrace $$ 
 
la base de l'espace vectoriel $\SS^{n}U$. On d\'efinit l'action de $\C^{*}$ sur $\SS^{n}U$ par
$$\left(\begin{array}{cc }  t^{\alpha}&0 \\0&    t^{\beta} \\ \end{array} \right) \in GL_{2}(\C)$$

o\`u $\alpha , \beta$ sont des entiers; cette action agit sur $ v_{p}$ comme suit
$$\left(\begin{array}{cc }  t^{\alpha}&0 \\0&    t^{\beta} \\ \end{array} \right).v_{p}=t^{n\alpha+p(\beta -\alpha)}. v_{p}.$$ 

Soit 
$$\BB:=\lbrace z_{p,q}:=x^{n-p}y^{p}\wedge x^{n-q}y^{q}  ,  \hspace{0.2 cm} 0\leq p< q\leq n \rbrace $$ 
 
la base de l'espace vectoriel $\bigwedge^{2}\SS^{n}U$. On d\'efinit l'action de $\C^{*}$ sur $\bigwedge^{2}\SS^{n}U$ par l'action pr\'ec\'edente de $\C^{*}$ qui agit sur $ z_{p,q}$ comme suit
$$\left(\begin{array}{cc }  t^{\alpha}&0 \\0&    t^{\beta} \\ \end{array} \right).z_{p,q}=t^{2n\alpha+(p+q)(\beta -\alpha)}. z_{p,q}.$$

Soient $k$ un entier avec $1\leq k\leq 2n-1$ et $E_{k}$ le sous-espace vectoriel de $\bigwedge^{2}\SS^{n}U $ tels que 
$$\left(\begin{array}{cc } t^{\alpha}&0 \\0&    t^{\beta} \\ \end{array} \right).u=t^{2n\alpha +k(\beta -\alpha)}. u,$$

pour tout $u\in E_{k}$. Alors on obtient
$$\bigwedge^{2}\SS^{n}U \simeq \bigoplus_{1\leq k\leq 2n-1}E_{k},$$

o\`u $E_{k}$ est engendr\'e par les \'el\'ements $z_{p,q}$ tels que $k=p+q$. 

\subsection{\bf Th\'eor\`eme} \label{3.4} 
{\em On utilise les m\^emes notations de la remarque \ref{3.3}. Soit $\P^{n}=\P(\SS^{n}U)$ l'espace projectif associ\'e \`a l'espace vectoriel $\SS^{n}U$. 

1- Soit $\WW$ l'ensemble de tous les sous-espaces vectoriels $W \subset \bigwedge^{2}(\SS^{n}U)$ tels que $W$ est \mbox{$\C^{*}$-invariant} et v\'erifie la condition (*). Alors il existe un ensemble $\ZZ_{k} \subset \P(E^{*}_{k})$ non vide pour tout entier $k$ avec $3 \leq k \leq 2n-3$ tel que 
$$\WW=\lbrace  \bigoplus_{3\leq k\leq 2n-3}W_{k}| W_{k}\in \ZZ_{k} ,  \hspace{0.2 cm} 3 \leq k \leq 2n-3 \rbrace. $$

2- Soit $W\in \WW$. Alors il existe un sous-espace vectoriel $D_{W} \subset \bigwedge^{2}\SS^{n}U$ tel que $D_{W}$ est \mbox{$\C^{*}$-invariant} et v\'erifie que
$$\bigwedge^{2}\SS^{n}U\simeq D_{W}\oplus W $$ 

est un $\C^{*}$-isomorphisme. De plus on a la suite exacte suivante
$$0\longrightarrow Q(-1)\stackrel{\varpi_{W}}{\longrightarrow} D_{W} \otimes\ko_{\P^{n}}  \longrightarrow  F(W)(1)   \longrightarrow 0,$$

o\`u $F(W)$ est le fibr\'e de Tango associ\'e au sous-espace $W$.}

\begin{proof}  
1- Soit 
$$\BB:=\lbrace z_{p,q}:=x^{n-p}y^{p}\wedge x^{n-q}y^{q}  ,  \hspace{0.2 cm} 0\leq p< q\leq n \rbrace $$
 
la base de l'espace vectoriel $\bigwedge^{2}\SS^{n}U$. Dans la remarque \ref{3.3}, on a vu
$$\bigwedge^{2}\SS^{n}U \simeq \bigoplus_{1\leq k\leq 2n-1}E_{k},$$

o\`u $E_{k}$ est engendr\'e par les \'el\'ements $z_{p,q}$ tels que $k=p+q$. On va d\'emontrer, pour tout $u\in E_{k}$, que $u$ est d\'ecomposable si et seulement s'il existe un entier $max(k-n,0)\leq p\leq [\frac{k-1}{2}]$ et  \mbox{$\varepsilon \in \C$} tels que $u=\varepsilon .z_{p,k-p} $. Notons $\lbrace z_{p,k-p}\wedge z_{q,k-q} \rbrace_{max(k-n,0)\leq p<q\leq [\frac{k-1}{2}]}$ la base de l'image de $E_{k}\times E_{k}$ dans $\bigwedge^{4}\SS^{n}U$ par le morphisme canonique suivant
$$\bigwedge^{2}\SS^{n}U\times \bigwedge^{2}\SS^{n}U\longrightarrow \bigwedge^{4}\SS^{n}U.$$

Si $u=\sum_{max(k-n,0)\leq p\leq [\frac{k-1}{2}]} a_{p}. z_{p,k-p} $, alors on obtient 
$$u\wedge u = 2 \sum_{max(k-n,0)\leq p<q\leq [\frac{k-1}{2}]} a_{p}.a_{q}. z_{p,k-p} \wedge z_{q,k-q}. $$

Si $u\wedge u=0$, on obtient que $a_{p}.a_{q}=0$ pour tout $max(k-n,0) \leq p<q\leq [\frac{k-1}{2}]$. Donc il existe $\varepsilon \in \C$ tel que $u=\varepsilon .z_{p,k-p} $. Il en d\'ecoule, pour tout $1\leq k\leq 2n-1$, que la dimension maximale d'un sous-espace vectoriel de $E_{k}$ ne contenant aucun \'el\'ement d\'ecomposable non nul est $dim(E_{k})- 1$ et que de tels sous-espaces existent. Plus pr\'ecis\'ement, pour tout $1\leq k\leq 2n-1$ et pour tout entier $p$ tel que 
$$max(k-n,0)\leq p\leq [\frac{k-1}{2}],$$
 
l'ensemble des hyperplans $H\subset E_{k}$ tels que $z_{p,k-p} \in H $ est un hyperplan de $\P(E^{*}_{k})$. Cet ensemble est le ferm\'e de Zariski de $\P(E^{*}_{k})$ suivant 
$$\kj (z_{p,k-p})=\lbrace H\subset E_{k} |z_{p,k-p} \in H  \rbrace \subset \P(E^{*}_{k}). $$

Soit $\DD(z_{p,k-p})=\P(E^{*}_{k})\diagdown \kj (z_{p,k-p})$. On d\'efinit
$$\ZZ_{k}= \bigcap_{p=max(k-n,0)}^{[\frac{k-1}{2}]}\DD(z_{p,k-p})$$

qui est l'ouvert constitu\'e de tous les sous-espaces vectoriels de $E_{k}$ de dimension $dim(E_{k})- 1$ ne contenant pas d'\'el\'ement d\'ecomposable non nul de $\bigwedge^{2}\SS^{n}U$. On a donc $$\ZZ_{1}=\ZZ_{2}=\ZZ_{2n-1}=\ZZ_{2n-2}=\emptyset .$$ 

Soit $W_{k}\in \ZZ_{k}$ pour tout $3\leq k\leq 2n-3$. On consid\`ere le sous-espace vectoriel de $\bigwedge^{2}\SS^{n}U$  
$$W= \bigoplus_{3\leq k\leq 2n-3}W_{k}.$$
 
Alors $W$ est $\C^{*}$-invariant et 
$$dim(W)=\left(\begin{array}{cc } n+1 \\2 \end{array} \right)-(2n-1).$$

Il reste \`a d\'emontrer que $W$ ne contient aucun \'el\'ement d\'ecomposable non nul de $\bigwedge^{2}\SS^{n}U$. On ~veut d\'emontrer que $w= 0$, tout en consid\'erant que $w\in W$ tel que $w \wedge w=0$. On a
$$w = \sum_{k=3}^{2n-3}w_{k}, \text{ o\`u }w_{k}\in W_{k}.$$

Donc on obtient
$$w \wedge w= 2 \sum _{i=3}^{2n-4} \left( \sum _{j=1}^{2n-3-i} w_{i} \wedge w_{i+j} \right)+  \sum _{i=3 }^{2n-3}   w_{i} \wedge w_{i}. $$
 
L'action de $\C^{*}$ sur un \'el\'ement $w_{i} \wedge w_{j} \in \bigwedge^{4}\SS^{n}U $ \'etant la multiplication par $t^{4n\alpha+2(i+j)(\beta -\alpha)}$, alors on regroupe les \'el\'ements dans $w \wedge w$ suivant l'action de $\C^{*}$  en mettant ensemble les \'el\'ements $w_{i} \wedge w_{j}$ ayant le m\^eme $i+j$. On obtient
$$w \wedge w=\sum_{d=6}^{4n-6} g_{d},$$

o\`u
$$g_{2m}=   w_{m} \wedge w_{m}+ 2 \sum _{j=1}^{m-3}  w_{m-j} \wedge w_{m+j},\hspace{0.2 cm} g_{2m+1}=  2 \sum _{j=0}^{m-3}    w_{m-j} \wedge w_{m+1+j}, $$
$$g_{2n}=   w_{n} \wedge w_{n}+ 2 \sum _{j=1}^{n-3}   w_{n-j} \wedge w_{n+j}, $$
$$g_{4n-2m}=   w_{2n-m} \wedge w_{2n-m}+ 2 \sum _{j=1}^{m-3}   w_{2n-m-j} \wedge w_{2n-m+j}, $$
$$g_{4n-2m-1}=  2 \sum _{j=0}^{m-3}  w_{2n-m-j} \wedge w_{2n-m-1+j}, $$

o\`u $3\leq m\leq n-1$. On a aussi, pour tout $t\in \C^{*}$,
$$\left(\begin{array}{cc } t^{\alpha}&0 \\0&    t^{\beta} \\ \end{array} \right).(w \wedge w)=\sum_{d=6}^{4n-6}\left( \left(\begin{array}{cc } t^{\alpha}&0 \\0&    t^{\beta} \\ \end{array} \right). g_{d}\right),$$
$$\left(\begin{array}{cc } t^{\alpha}&0 \\0&    t^{\beta} \\ \end{array} \right).(w \wedge w)=\sum_{d=6}^{4n-6} (t^{4n\alpha+2d(\beta -\alpha)}. g_{d}).$$

Comme $w \wedge w=0$, alors on a $g_{d}=0$ pour tout $6\leq d \leq 4n-6$. Nous allons montrer, par la r\'ecurrence sur $m$, que $w_{r}=0,w_{2n-r}=0 $ pour tout $3\leq r \leq m$ et pour tout $3\leq m\leq n$:

- Pour $m=3$, on a $g_{6}=w_{3} \wedge w_{3}=0$ et $g_{4n-6}=w_{2n-3} \wedge w_{2n-3}=0$, ce qui entra\^ine que $w_{3}=w_{2n-3}=0$. 

- Pour $m=4$, on a $g_{6}=w_{3} \wedge w_{3}=0$ et $g_{4n-6}=w_{2n-3} \wedge w_{2n-3}=0$, ce qui entra\^ine que $w_{3}=w_{2n-3}=0$, et on obtient que $g_{8}=w_{4} \wedge w_{4}=0$ et $g_{4n-8}=w_{2n-4} \wedge w_{2n-4}=0$, ce qui entra\^ine que $w_{4}=w_{2n-4}=0$. 

- On suppose que $w_{r}=0,w_{2n-r}=0 $ pour tout $3\leq r \leq m$, et on montre que $w_{r}=0,w_{2n-r}=0 $ pour tout $3\leq r \leq m+1$. 

- Pour $m+1$, on a 
$$g_{2(m+1)}=w_{m+1} \wedge w_{m+1}+2  ( w_{m} \wedge w_{m+2}+$$
$$w_{m-1} \wedge w_{m+3}+ \ldots   +w_{4} \wedge w_{2m-2}+w_{3} \wedge w_{2m-1})=0,$$

 et 
$$g_{4n-2(m+1)}=w_{2n-(m+1)} \wedge w_{2n-(m+1)}+ 2  (  w_{2n-m-2} \wedge w_{2n-m}+$$ 
$$  w_{2n-m-3} \wedge w_{2n-m+1}+ \ldots +w_{2n-2m+2} \wedge w_{2n-4}+w_{2n-2m+1} \wedge w_{2n-3})=0,$$

ce qui entra\^ine que $w_{m+1}=w_{2n-(m+1)}=0$. Pour tout $3\leq r \leq m+1$, on a alors $w_{r}=w_{2n-r}=0 $. Donc on obtient que $w= 0$, et que $W$ ne contient aucun \'el\'ement d\'ecomposable non nul de $\bigwedge^{2}\SS^{n}U$. On en d\'eduit
$$\lbrace  \bigoplus_{3\leq k\leq 2n-3}W_{k}| W_{k}\in \ZZ_{k} ,  \hspace{0.2 cm} 3 \leq k \leq 2n-3 \rbrace \subseteq \WW. $$

On d\'emontre maintenant l'inclusion inverse. Soient un entier $k$ tel que $1\leq k\leq 2n-1$, et $W\in \WW$ qui est $\C^{*}$-invariant. On d\'efinit les espaces $W_{k}$ par $u_{k} \in E_{k} $ tels que $\sum_{1\leq k\leq 2n-1}u_{k} \in W$. On en d\'eduit
$$W = \bigoplus_{1\leq k\leq 2n-1}W_{k}.$$

Comme le sous-espace vectoriel $E_{1}$ (resp. $E_{2},E_{2n-1},E_{2n-2}$) contient un seul \'el\'ement $z_{p,q}$ tel que $p+q=1$ (resp. $p+q=2,2n-1,2n-2$) et comme $W$ ne contient aucun \'el\'ement d\'ecomposable non nul de $\bigwedge^{2}\SS^{n}U$, alors $W_{1}=W_{2}=W_{2n-1}=W_{2n-2}=0$. On obtient donc que 
$$W = \bigoplus_{3\leq k\leq 2n-3}W_{k}.$$

Comme $W$ ne contient aucun \'el\'ement d\'ecomposable non nul de $\bigwedge^{2}\SS^{n}U$, alors $W_{k}$ ne contient aucun \'el\'ement d\'ecomposable non nul de $\bigwedge^{2}\SS^{n}U$. Donc $W_{k}\in \ZZ_{k}$ pour tout $3\leq k\leq 2n-3$ et on a
$$\WW\subseteq \lbrace  \bigoplus_{3\leq k\leq 2n-3}W_{k}| W_{k}\in \ZZ_{k} ,  \hspace{0.2 cm} 3 \leq k \leq 2n-3 \rbrace . $$ 

2- Il suffit de choisir $d_{k}\in E_{k}\diagdown W_{k}$, pour tout $1\leq k\leq 2n-1$, et de consid\'erer le sous-espace vectoriel de $\bigwedge^{2}\SS^{n}U$  

$$D_{W}= \bigoplus_{1\leq k\leq 2n-1}\C.d_{k}.$$

D'apr\`es la d\'efinition \ref{3.1}, on obtient la suite exacte recherch\'ee.

\end{proof}

\subsection{\bf Proposition} {\em (D\'ecomposition Clebsch-Gordan).\label{3.5}
Soit $U$ un espace vectoriel complexe de dimension $2$. Alors il existe une $SL_{2}(\C)$-d\'ecomposition irr\'eductible de $
\SS^{n}U\otimes \SS^{n}U$
$$\SS^{n}U\otimes \SS^{n}U \simeq  \bigoplus_{i=0}^{n} \SS^{2n-2i}U .$$

En particulier, on a
$$\bigwedge^{2}(\SS^{n}U)\simeq \SS^{2(n-1 )}U \oplus \SS^{2(n- 3)}U \oplus \SS^{2(n- 5)}U \oplus \SS^{2(n- 7)}U \oplus \ldots $$

qui est $SL_{2}(\C)$-isomorphisme.}

\begin{proof}
Voir \cite{kr-pr} pages 93-97. 

\end{proof}

\subsection{\bf Proposition} \label{3.6} 
 {\em On utilise les m\^emes notations de la remarque \ref{3.3}. Soient $\beta=-\alpha$ et le sous-espace vectoriel 
$$W=\SS^{2(n- 3)}U \oplus \SS^{2(n- 5)}U \oplus \SS^{2(n- 7)}U \oplus \ldots \subset \bigwedge^{2}(\SS^{n}U). $$

Alors $W$ est $SL_{2}(\C)$-invariant et $W\in \WW$. Le fibr\'e de Tango $F(W)$, qui est d\'efini par la suite exacte suivante
$$0\longrightarrow Q(-1)\stackrel{\varpi_{W}}{\longrightarrow} \SS^{2(n-1 )}U \otimes\ko_{\P^{n}}  \longrightarrow  F(W)(1)   \longrightarrow 0,$$

est $SL_{2}(\C)$-invariant.}

\begin{proof}
Voir l'article de Cascini \cite{ca}, proposition 2.1. 

\end{proof}

\subsection{\bf Proposition}\label{3.7}
{\em On utilise les m\^emes notations de la remarque \ref{3.3}. Soient $i,n, \alpha ,\gamma \in \N$ et $\beta \in \Z$ tels que $n>2$, $\gamma >0$, $\alpha \geq \beta$, $\alpha +\beta\geq 0$ et $\gamma + n\alpha+i(\beta -\alpha)> 0$ pour $0\leq i\leq n$. Soient $g_{0}, \ldots ,g_{n}$ des formes homog\`enes sans z\'ero commun sur $\P(\SS^{n}U)$ telles que 
$$ deg(g_{i})=\gamma +n\alpha+ i(\beta -\alpha),\hspace{0.2 cm} i=0,1, \ldots ,n.$$
 
Soient $Q$ le fibr\'e de quotient et $F(W)$ le fibr\'e de Tango sur $\P(\SS^{n}U)$ pour $  W  \in \WW$ et $D_{W}$ comme dans le th\'eor\`eme \ref{3.4}. Les fibr\'es $Q$ et $F(W)$ sont d\'efinis par les suites exactes
$$ 0\longrightarrow \ko_{\P(\SS^{n}U)}(-1)\stackrel{g}{\longrightarrow} \SS^{n}U\otimes \ko_{\P(\SS^{n}U)} \longrightarrow Q \longrightarrow 0,$$ 

$$0\longrightarrow Q(-1)\stackrel{\varpi_{W}}{\longrightarrow} D_{W} \otimes\ko_{\P(\SS^{n}U)}  \longrightarrow  F(W)(1)   \longrightarrow 0,$$

o\`u $g$ est le morphisme canonique. Alors le fibr\'e $Q$ (resp. $F(W)$) poss\`ede une image invers\'ee g\'en\'eralis\'ee $Q_{\gamma,\alpha,\beta}$ (resp. $F_{\gamma,\alpha,\beta}$) d\'efinie par

$$  0\longrightarrow \ko_{\P^{n}}(-\gamma)\stackrel{g(-\gamma)}{\longrightarrow} \SS^{n}\ku\longrightarrow \kq_{\gamma,\alpha,\beta} \longrightarrow 0$$

$$\text{ (resp. } 0\longrightarrow \kq_{\gamma,\alpha,\beta}(-\gamma) \stackrel{\varpi_{W}(-2\gamma)}{\longrightarrow}\SS^{2(n-1)}\ku \longrightarrow  \kf_{\gamma,\alpha,\beta}(\gamma) \longrightarrow 0),$$

o\`u $\ku =\ko_{\P^{n}}(\alpha )\oplus \ko_{\P^{n}}(\beta)$, et $\kf_{\gamma,\alpha,\beta}(\gamma):=F_{\gamma,\alpha,\beta}(-2\gamma)$, $\kq_{\gamma,\alpha,\beta}:=Q_{\gamma,\alpha,\beta}(-\gamma)$. La premi\`ere classe de Chern du fibr\'e $\kf_{\gamma,\alpha,\beta}$ est $n(\beta +\alpha)(2n-1-(\dfrac{n+1}{2}))$. On appelle le~fibr\'e $\kq_{\gamma,\alpha,\beta}$ le fibr\'e de quotient 
pond\'er\'e par les poids $\gamma,\alpha,\beta$, provenant d'une image inverse g\'en\'eralis\'ee sur $\P^{n}$. On appelle le fibr\'e $\kf_{\gamma,\alpha,\beta} $ le fibr\'e de Tango pond\'er\'e par les poids $\gamma,\alpha,\beta$, provenant d'une image inverse g\'en\'eralis\'ee sur $\P^{n}$}.

\begin{proof} On consid\`ere l'application 
$$\omega:= (g_{0},\ldots ,g_{n}): \C^{n+1}\setminus  \lbrace 0 \rbrace \longrightarrow \SS^{n}U \setminus  \lbrace 0 \rbrace \hspace{4 cm}$$
$$\hspace{2 cm}v\longmapsto (g_{0}(v),\ldots ,g_{n}(v)).$$

On consid\`ere l'action de $\C^{*}$ sur $\SS^{n}U $
$$\sigma:\C^{*}\times  \SS^{n} U \longrightarrow \SS^{n} U   $$
$$(t,u)\longmapsto t^{\gamma}.\left(\begin{array}{cc } t^{\alpha}&0 \\0&    t^{\beta} \\ \end{array} \right).u$$

qui est repr\'esent\'ee par la matrice
$$ t^{\gamma}. \left(
\begin{array}{ccccccccc}
 t^{n\alpha} &  &   & & &            &   \\
     & t^{n\alpha+(\beta -\alpha)}  &  &   & &   &\bigzero \\
  &  & \ddots &   &    &   &   \\
\bigzero &  &   &    &  t^{n\alpha+(n-1)(\beta -\alpha)} &  &   \\
 &  &   &    &    &  &      t^{n\beta} \\
\end{array}
\right)\in PGL(\SS^{n} U).
$$

On consid\`ere aussi l'action de $\C^{*}$ sur $\C^{n+1} $ qui est la multiplication usuelle sur $\C^{n+1} $
$$T:\C^{*}\times  \C^{n+1} \longrightarrow \C^{n+1}   $$
$$(t,u)\longmapsto t.u.$$

Alors $\omega $ est une $\C^{*}$-application par rapport \`a ces deux actions. L'action $\sigma$ induit une action $\overline{ \sigma}\in PGL(\SS^{n} U) $ de $\C^{*}$ sur $\P(\SS^{n} U)$ et l'action $T$ induit une action triviale de $\C^{*}$ sur $\P^{n} $. Donc la transform\'ee de Horrocks \ref{2.2} est
\xmat{{\bf Iminvg}:  \kf\kv(\P(\SS^{n} U),\overline{\sigma})   \ar[rr] &  &  \kf\kv(\P^{n}),  \\
 }

o\`u $\kf\kv(\P(\SS^{n} U),\overline{\sigma})$ est la cat\'egorie de fibr\'es vectoriels sur $\P(\SS^{n} U)$ qui sont $\C^{*}$-invariants \mbox{au-dessus} de l'action $\overline{ \sigma}$ (resp. $ \kf\kv(\P^{n})$ est la cat\'egorie de fibr\'es vectoriels sur $\P(\SS^{n} U)$ qui sont $\C^{*}$-invariants au-dessus de l'action $T$). On consid\`ere le morphisme
$$g:={}^{T}\omega=:\ko_{\P(\SS^{n}U)}(-1) \longrightarrow \SS^{n}U \otimes \ko_{\P(\SS^{n}U)} ,$$

avec $\ko_{\P(\SS^{n}U)}(-1)$ et $\SS^{n}U \otimes \ko_{\P(\SS^{n}U)}$ qui sont munis de l'action canonique $\sigma(t)$. Alors $g$ est un $\C^{*}$-morphisme. Comme on a, pour tout $t\in \C^{*}$,
$$\sigma(t).(x^{n-i}.y^{i}) =t^{n\alpha+i(\beta -\alpha) +\gamma}.(x^{n-i}.y^{i}) $$
 
alors le sous-fibr\'e $(x^{n-i}.y^{i}.\C)\otimes \ko_{\P(\SS^{n} U)} $ de $\SS^{n}U \otimes \ko_{\P(\SS^{n}U)}$ est $\C^{*}$-invariant. On obtient
$$\ko_{\P(\SS^{n} U)}^{(n\alpha+i(\beta -\alpha) +\gamma)}\simeq(x^{n-i}.y^{i}.\C)\otimes \ko_{\P(\SS^{n} U)}  $$

qui est d\'efini localement, pour tout $v\in \SS^{n} U $, par
$$(\ko_{\P(\SS^{n} U)}^{(n\alpha+i(\beta -\alpha) +\gamma)})_{v} \simeq \C \stackrel{\simeq}{\longrightarrow}((x^{n-i}.y^{i}.\C) \otimes \ko_{\P(\SS^{n} U)})_{v}\simeq x^{n-i}.y^{i}.\C   $$
$$a\longmapsto a x^{n-i}.y^{i}.\hspace{1.5 cm}$$

Donc on a un $\C^{*}$-isomorphisme
$$\SS^{n}U \otimes \ko_{\P(\SS^{n}U)}\simeq \bigoplus_{i=0}^{n}\ko_{\P(\SS^{n} U)}^{(n\alpha+i(\beta -\alpha) +\gamma)}. $$

Alors on a un $\C^{*}$-morphisme
$$g:={}^{T}\omega=:\ko_{\P(\SS^{n}U)}(-1) \longrightarrow \bigoplus_{i=0}^{n}\ko_{\P(\SS^{n} U)}^{(n\alpha+i(\beta -\alpha) +\gamma)} .$$

On consid\`ere le morphisme 
$$\varpi_{W}:Q (-1) \stackrel{ }{\longrightarrow} D_{W}\otimes \ko_{\P(\SS^{n}U)},$$ 

avec $D_{W}\otimes \ko_{\P(\SS^{n}U)}$ et $Q(-1)$ qui sont munis de l'action canonique $\sigma(t)$. Pour que le morphisme $\varpi_{W}$ soit un $\C^{*}$-morphisme il faut avoir
$$(\varpi_{W})_{v}(\sigma(t).v_{1})=t^{q}\sigma(t).(\varpi_{W})_{v}(v_{1})$$

pour tout $v,v_{1}\in \SS^{n}U$, $q$ un entier et
$$(\varpi_{W})_{v}:Q (-1)_{v}  \stackrel{ }{\longrightarrow} (D_{W}\otimes \ko_{\P(\SS^{n}U)})_{v}=D_{W}= \bigoplus_{1\leq k\leq 2n-1}\C.d_{k},$$

o\`u $d_{k}\in E_{k}\diagdown W_{k}$, pour tout $1\leq k\leq 2n-1$, comme dans le th\'eor\`eme \ref{3.4}. \mbox{Alors on obtient $q=0$}. Comme on a, pour tout $t\in \C^{*}$,
$$\sigma(t).(d_{k}) =t^{2n\alpha+k(\beta -\alpha) +2\gamma}.(d_{k}), $$
 
alors le sous-fibr\'e $(d_{k}.\C)\otimes \ko_{\P(\SS^{n} U)} $ de $D_{W} \otimes \ko_{\P(\SS^{n}U)}$ est $\C^{*}$-invariant. On a
$$\ko_{\P(\SS^{n} U)}^{(2n\alpha+k(\beta -\alpha) +2\gamma)}\simeq (d_{k}.\C)\otimes \ko_{\P(\SS^{n} U)} $$

qui est d\'efini localement, pour tout $v\in \SS^{n} U $, par
$$(\ko_{\P(\SS^{n} U)}^{((2n\alpha+k(\beta -\alpha) +2\gamma)})_{v} \simeq \C \stackrel{\simeq}{\longrightarrow}((d_{k}.\C)\otimes \ko_{\P(\SS^{n} U)})_{v}\simeq d_{k}.\C  . $$
$$a\longmapsto a d_{k}.\hspace{0.1 cm}$$

Donc on en d\'eduit le $\C^{*}$-isomorphisme suivant
$$D_{W} \otimes \ko_{\P(\SS^{n}U)} \simeq \bigoplus_{k=1}^{2n-1}\ko_{\P(\SS^{n} U)}^{(2n\alpha+k(\beta -\alpha) +2\gamma)}, $$

et le $\C^{*}$-morphisme
$$\varpi_{W}:Q(-1) \stackrel{ }{\longrightarrow} \bigoplus_{k=1}^{2n-1}\ko_{\P(\SS^{n} U)}^{(2n\alpha+k(\beta -\alpha)+2\gamma)}.$$

Comme le fibr\'e $Q$ a une $GL(\SS^{n} U)$-action alors il a une $\C^{*}$-action. Autrement dit, le fibr\'e $Q$ est \mbox{$\C^{*}$-invariant} au-dessus de l'action $\sigma(t)$. D'apr\`es la proposition \ref{3.2}, on obtient que
$$\sigma(t)^{*}F(W) \simeq F(W). $$ 

D'apr\`es la d\'efinition \ref{2.2} et la proposition \ref{2.3}, on obtient 
$${\bf Iminvg}(Q(-1))=Q_{\gamma,\alpha,\beta} \text{ et } {\bf Iminvg}(F(W)(1))=F_{\gamma,\alpha,\beta} $$ 

et on a aussi
$${\bf Iminvg}(\ko_{\P(\SS^{n}U)}(-1))= \ko_{\P^{n}},$$

$${\bf Iminvg}(D_{W} \otimes \ko_{\P(\SS^{n}U)} )\simeq {\bf Iminvg}(\bigoplus_{k=1}^{2n-1}\ko_{\P(\SS^{n} U)}^{(2n\alpha+k(\beta -\alpha) +2\gamma)})$$
$$=\bigoplus_{k=1}^{2n-1}\ko_{\P^{n}}(2n\alpha+k(\beta -\alpha)+2\gamma):=\SS^{2(n-1)}\ku (2\gamma) ,$$

$${\bf Iminvg}(\SS^{n}U \otimes \ko_{\P(\SS^{n}U)})\simeq {\bf Iminvg}( \bigoplus_{i=0}^{n}\ko_{\P(\SS^{n} U)}^{(n\alpha+i(\beta -\alpha) +\gamma)})$$
$$=\bigoplus_{i=0}^{n}\ko_{\P^{n}}(n\alpha+i(\beta -\alpha)+\gamma):= \SS^{n}\ku (\gamma) ,$$

o\`u $\ku =\ko_{\P^{n}}(\alpha )\oplus \ko_{\P^{n}}(\beta)$. On en d\'eduit les suites exactes
$$ 0\longrightarrow \ko_{\P^{n}}\stackrel{g}{\longrightarrow} \bigoplus_{i=0}^{n}\ko_{\P^{n}}(n\alpha+i(\beta -\alpha)+\gamma) \longrightarrow Q_{\gamma,\alpha,\beta} \longrightarrow 0$$

et  
$$0\longrightarrow Q_{\gamma,\alpha,\beta} \stackrel{\varpi_{W}}{\longrightarrow}\bigoplus_{k=1}^{2n-1}\ko_{\P^{n}}(2n\alpha+k(\beta -\alpha)+2\gamma) \longrightarrow  F_{\gamma,\alpha,\beta} \longrightarrow 0.$$

Ces suites exactes s'\'ecrivent \'egalement comme suit
$$ 0\longrightarrow \ko_{\P^{n}}(-\gamma)\stackrel{g(-\gamma)}{\longrightarrow} \bigoplus_{i=0}^{n}\ko_{\P^{n}}(n\alpha+i(\beta -\alpha)) \longrightarrow \kq_{\gamma,\alpha,\beta} \longrightarrow 0$$

et  
$$0\longrightarrow \kq_{\gamma,\alpha,\beta}(-\gamma) \stackrel{\varpi_{W}(-2\gamma)}{\longrightarrow}\bigoplus_{k=1}^{2n-1}\ko_{\P^{n}}(2n\alpha+k(\beta -\alpha)) \longrightarrow  \kf_{\gamma,\alpha,\beta}(\gamma) \longrightarrow 0,$$

o\`u $\kq_{\gamma,\alpha,\beta}:=Q_{\gamma,\alpha,\beta}(-\gamma)$, et $\kf_{\gamma,\alpha,\beta}(\gamma):=F_{\gamma,\alpha,\beta}(-2\gamma)$.

\end{proof}

\subsection{\bf Remarque}\label{3.8} 
Le fibr\'e $\kf_{\gamma,\alpha,-\alpha}$ pond\'er\'e par les poids $\gamma,\alpha ,-\alpha$ est le fibr\'e pond\'er\'e par les poids $\gamma,\alpha$ de Cascini \cite{ca}. On va noter $\kq=\kq_{\gamma,\alpha,\beta}$ et $\kf=\kf_{\gamma,\alpha,\beta}$.

\newpage

\section{ \large\bf Stabilit\'e et d\'eformation miniversale du fibr\'e de Tango pond\'er\'e.}
 
\vspace{1cm}

Nous allons d\'emontrer que le fibr\'e $\kf$ est stable et que les fibr\'es $\kq$ et $\kf$ sont invariants par rapport \`a une d\'eformation miniversale. Nous allons aussi montrer que l'espace de Kuranishi du fibr\'e $\kf$ est lisse au point correspondant du fibr\'e $\kf$.

\subsection{\bf Proposition}\label{3.9} 
{\em Soit $\kf$ le fibr\'e de Tango pond\'er\'e par les poids $\gamma,\alpha,\beta$ sur $\P^{n}$. 

1- Si on a $\gamma > 2n\alpha+(\beta-\alpha)$, alors $\kf$ est stable.

2- Soit $\gamma> n\alpha$. Si $\kf$ est stable, alors on a $\gamma > 2n\alpha+(\beta-\alpha)$.}
\begin{proof}
Soient $q$ un entier avec $1\leq q\leq n-2$ et $t\in \Z$ tels que $(\bigwedge^{q}\kf)_{norm}=\bigwedge^{q}\kf (t)$. On obtient alors
$$\dfrac{c_{1}(\bigwedge^{q}\kf (t))}{rg(\bigwedge^{q}\kf (t))}\leq 0,$$
qui s'\'ecrit \'egalement
$$t+c_{1}(\kf)\dfrac{q}{n-1}\leq 0.$$

Donc on a
$$t\leq -\dfrac{n}{n-1}(\alpha +\beta)(2n-1-(\frac{n+1}{2}))< (\alpha +\beta)(2-n)\leq 0.$$

1- De la suite exacte suivante
$$  0\longrightarrow \ko_{\P^{n}}(-\gamma)\stackrel{}{\longrightarrow} \SS^{n}\ku\longrightarrow \kq \longrightarrow 0,$$

il en d\'ecoule la suite exacte suivante
$$0\longrightarrow \SS^{q-1} \SS^{n}\ku(-\gamma+m) \longrightarrow \SS^{q} \SS^{n}\ku(m)\longrightarrow \SS^{q}\kq(m) \longrightarrow 0$$

qui nous donne $h^{i}(\SS^{q}\kq(m))=0 $, pour tout $m\in \Z$ et $1\leq i\leq n-2$. De la suite exacte suivante
$$0\longrightarrow \kq(-\gamma) \stackrel{}{\longrightarrow}\SS^{2(n-1)}\ku \longrightarrow  \kf(\gamma) \longrightarrow 0,$$

on obtient la r\'esolution suivante
$$0\longrightarrow \SS^{q} \kq (-2q\gamma+t) \longrightarrow \SS^{q-1} \kq \otimes \SS^{2(n-1)}\ku(-\gamma(2q-1) +t) 
$$
$$\stackrel{}{\longrightarrow} \SS^{q-2} \kq\otimes \bigwedge^{2} \SS^{2(n-1)}\ku(-\gamma(2q-2) +t) \stackrel{a_{q-2}}{\longrightarrow} \ldots$$
$$\hspace{2 cm} \ldots 
\stackrel{a_{2}}{\longrightarrow}   \kq \otimes \bigwedge^{q-1} \SS^{2(n-1)}\ku(-\gamma(q+1) +t) \stackrel{a_{1}}{\longrightarrow}  \bigwedge^{q} \SS^{2(n-1)}\ku(-q\gamma +t) 
\stackrel{a_{0}}{\longrightarrow} \bigwedge^{q} \kf(t) \longrightarrow 0.$$

En consid\'erant $A_{j}=ker (a_{j})$ pour $j=0,1, \ldots ,q-2$, on obtient alors $h^{i}(A_{0})=0$ pour tout~\mbox{ $1\leq i\leq n+1-q$}. On a
$$q\gamma -t>2nq\alpha +q(\beta-\alpha)+(n-2)(\beta+\alpha)\geq 2nq\alpha +q(\beta-\alpha)  $$
$$\geq max\lbrace e \in \Z|\hspace{0.2 cm}  \ko_{\P^{2n+1}}(e)\subseteq \bigwedge^{q} \SS^{2(n-1)}\ku  \rbrace = 2n q\alpha +(\beta-\alpha)\dfrac{q(q+1)}{2}.$$

De la suite exacte suivante
$$0\longrightarrow A_{0}\longrightarrow  \bigwedge^{q} \SS^{2(n-1)}\ku(-q\gamma +t) \stackrel{a_{0}}{\longrightarrow} \bigwedge^{q} \kf(t) \longrightarrow 0,$$

on en d\'eduit $ h^{0}(\bigwedge^{q} \kf(t))=0$. D'apr\`es le crit\`ere de Hoppe \cite{hopp}, $\kf$ est stable.

2- Supposons que $\gamma \leq 2n\alpha+(\beta-\alpha)$. On obtient 
$$ \gamma -t < 2n\alpha+(\beta-\alpha),$$

et on a aussi 
$$2\gamma -t > 2n\alpha+(n-2)(\beta+\alpha)\geq n\alpha.$$

Des suites exactes suivantes
$$  0\longrightarrow \ko_{\P^{n}}(-3\gamma+t)\stackrel{}{\longrightarrow} \SS^{n}\ku(-2\gamma+t)\longrightarrow \kq (-2\gamma+t) \longrightarrow 0$$

$$0\longrightarrow \kq(-2\gamma+t) \stackrel{}{\longrightarrow}\SS^{2(n-1)}\ku (-\gamma+t) \longrightarrow  \kf(t) \longrightarrow 0,$$ 

on en d\'eduit que $h^{0}(\kq (-2\gamma+t))=0$ et $h^{0}(\kf (t))=h^{0}(\SS^{2(n-1)}\ku (-\gamma+t))\neq 0$. Donc il existe un fibr\'e en droite trivial dans $\kf (t)$
$$\ko_{\P^{n}}\hookrightarrow \kf (t),$$

ce qui nous donne $\dfrac{c_{1}(\ko_{\P^{n}})}{rg(\ko_{\P^{n}})}=0\geq \dfrac{c_{1}(\kf (t))}{rg(\kf (t))}$. Donc $\kf $ n'est pas stable.

\end{proof}

\subsection{\bf Th\'eor\`eme}\label{4.1}
{\em (Hartshorne, \cite{ha10}). Si $E$ est un faisceau coh\'erent sur un sch\'ema projectif~$X$ sur un corps de base $K$ tel que $hdE\leq 1$, il existe un sch\'ema $Y=Spec(R)$ qui param\'etrise les~d\'eformations miniversales de $E$, o\`u $R$ est une $K$-alg\`ebre locale compl\`ete.}
 
\begin{proof}
Voir le th\'eor\`eme (19.1 \cite{ha10}).
\end{proof}

\subsection{\bf Th\'eor\`eme} 
{\em Soit $E$ un fibr\'e vectoriel sur la vari\'et\'e alg\'ebrique $\P^{n}$. Il existe un espace de Kuranishi $Kur(E)$ de $E$ qui est une base de la d\'eformation miniversale de $E$ ($Kur(E)$ param\'etrise toutes les d\'eformations miniversales de $E$).}

\begin{proof}
Voir l'article de M. Kuranishi \cite{kur}.
\end{proof}

Soit $e$ un point correspondant au fibr\'e $E$. Alors l'espace $Kur(E)$ est \'equip\'e d'une famille \mbox{universelle} et la fibre $(Kur(E),e)$, un espace topologique point\'e, est unique \`a un automorphisme pr\`es.

\subsection{\bf Lemme}\label{4.3} 
{\em Soient  $ \kq^{'}$ et $ \kq^{''}$ deux fibr\'es de quotient pond\'er\'es par les poids $\gamma,\alpha,\beta$ sur $\P^{n}$ qui sont d\'efinis par les suites exactes suivantes

$$0\longrightarrow \ko_{\P^{n}}(-\gamma)\longrightarrow \SS^{n}\ku\stackrel{q_{1}}{\longrightarrow} \kq^{'} \longrightarrow 0$$

$$0\longrightarrow \ko_{\P^{n}}(-\gamma)\longrightarrow \SS^{n}\ku\stackrel{q_{2}}{\longrightarrow} \kq^{''} \longrightarrow 0,$$

tels qu'il existe un morphisme $\psi : \kq^{'} \longrightarrow \kq^{''}$. Alors 
il existe un morphisme $\varphi : \SS^{n}\ku  \longrightarrow \SS^{n}\ku $ tel que 
$q_{2}\circ \varphi= \psi\circ q_{1}$}.

\begin{proof}
La d\'emonstration de ce lemme est tr\`es similaire \`a celle du lemme 4.3 \cite{bah}. 
\end{proof}

\subsection{\bf Lemme}\label{4.4} 
{\em Soient  $f,\hspace{0.2cm}f^{'}\in Hom(\ko_{\P^{n}}(-\gamma),\SS^{n}\ku) $ deux 
morphismes. Alors $f$ et \mbox{$f'$ donnent} le m\^eme \'el\'ement dans le sch\'ema 
$\QQuot_{\SS^{n}\ku/\P^{n}}$ si et seulement s'il existe un isomorphisme \mbox{$g \in End( 
\ko_{\P^{n}}(-\gamma))$} tel que $f=f^{'}\circ g$}.

\begin{proof}
C'est la d\'efinition du sch\'ema $\QQuot_{\SS^{n}\ku/\P^{n}}$.
  
\end{proof}

\subsection{\bf Th\'eor\`eme }\label{4.5} 
{\em Soit $ \kq_{0}$ un fibr\'e de quotient pond\'er\'e par les poids $\gamma,\alpha,\beta$ sur $\P^{n}$ qui est d\'efini par la suite exacte suivante
$$0\longrightarrow \ko_{\P^{n}}(-\gamma)\stackrel{x_{0}}{\longrightarrow} \SS^{n}\ku\longrightarrow \kq_{0} \longrightarrow 0$$

o\`u $x_{0}\in Hom(\ko_{\P^{n}}(-\gamma),\SS^{n}\ku)$.
Alors chaque d\'eformation miniversale du fibr\'e $\kq_{0}$ est encore un fibr\'e de quotient pond\'er\'e sur $\P^{n}$. L'espace de Kuranishi de $\kq_{0}$ est lisse au point correspondant de $\kq_{0}$.}

\begin{proof}
La d\'emonstration de ce th\'eor\`eme est tr\`es similaire \`a celle du th\'eor\`eme 4.5 ~\cite{bah}.
\end{proof}

\subsection{\bf Lemme}\label{4.6} 
{\em Soit $\kq$ un fibr\'e de quotient pond\'er\'e par les poids $\gamma,\alpha,\beta$ sur $\P^{n}$. Soient $\kf^{'}$ et $ \kf^{''}$ des fibr\'es de Tango pond\'er\'es par les poids $\gamma,\alpha,\beta$ sur $\P^{n}$ qui sont d\'efinis par les suites exactes
$$0\longrightarrow \kq(-\gamma) \stackrel{}{\longrightarrow}\SS^{2(n-1)}\ku \stackrel{p_{1}}{\longrightarrow}  \kf^{'}(\gamma) \longrightarrow 0$$

et
$$0\longrightarrow \kq(-\gamma) \stackrel{}{\longrightarrow}\SS^{2(n-1)}\ku \stackrel{p_{2}}{\longrightarrow} \kf^{''}(\gamma) \longrightarrow 0,$$

tels qu'il existe un morphisme $\psi :\kf^{'} (\gamma)\longrightarrow \kf^{''} (\gamma)
$. Alors il existe un morphisme 
$$\varphi :\SS^{2(n-1)}\ku  \longrightarrow \SS^{2(n-1)}\ku   ,$$ 

tel que $p_{2}\circ \varphi= \psi\circ p_{1}$.}

\begin{proof}
 De la suite exacte suivante
$$0\longrightarrow \kq(-\gamma) \stackrel{}{\longrightarrow}\SS^{2(n-1)}\ku \stackrel{p_{2}}{\longrightarrow} \kf^{''}(\gamma) \longrightarrow 0,$$

il en d\'ecoule la suite exacte suivante de groupes cohomologiques
$$0\longrightarrow Hom(\SS^{2(n-1)}\ku,\kq(-\gamma)) \longrightarrow 
Hom(\SS^{2(n-1)}\ku,\SS^{2(n-1)}\ku)\stackrel{p_{2}\hspace{0.1 cm}\circ \hspace{0.1 
cm} \bullet }{\longrightarrow} Hom(\SS^{2(n-1)}\ku,\kf^{''}(\gamma)) $$
$$\longrightarrow Ext^{1}(\SS^{2(n-1)}\ku,\kq(-\gamma)) 
\longrightarrow 0.$$
 
De la suite exacte suivante

$$  0\longrightarrow \ko_{\P^{n}}(-2\gamma)\otimes (\SS^{2(n-1)}\ku)^{*} \stackrel{}{\longrightarrow} \SS^{n}\ku(-\gamma)\otimes (\SS^{2(n-1)}\ku)^{*}\longrightarrow \kq(-\gamma) \otimes (\SS^{2(n-1)}\ku)^{*} \longrightarrow 0$$

on en d\'eduit

$$Ext^{1}(\SS^{2(n-1)}\ku,\kq(-\gamma))=H^{1}(\kq(-\gamma) \otimes (\SS^{2(n-1)}\ku)^{*} )=0.$$

Alors on obtient que 
$$Hom(\SS^{2(n-1)}\ku,\SS^{2(n-1)}\ku)\stackrel{p_{2}\hspace{0.1 cm}\circ \hspace{0.1 
cm} \bullet }{\longrightarrow} Hom(\SS^{2(n-1)}\ku,\kf^{''}(\gamma))\longrightarrow 0.$$

Comme le morphisme $\psi\circ p_{1}:\SS^{2(n-1)}\ku\longrightarrow 
\kf^{''}(\gamma) $ appartient \`a $ Hom(\SS^{2(n-1)}\ku,\kf^{''}(\gamma))$, alors il existe un 
morphisme $\varphi: \SS^{2(n-1)}\ku\longrightarrow \SS^{2(n-1)}\ku$ tel que 
$p_{2}\circ \varphi= \psi\circ p_{1}$.
 
\end{proof}

\subsection{\bf Lemme} \label{4.7}
{\em Soit $\kq$ un fibr\'e de quotient pond\'er\'e sur $\P^{n}$. Soient $f$ et $ f^{'}$ deux morphismes dans $Hom(\kq(-\gamma),\SS^{2(n-1)}\ku) $. Alors les morphismes $f$ et $f^{'}$ donnent le m\^eme \'el\'ement dans $\QQuot_{\SS^{2(n-1)}\ku/\P^{n}}$ si~et seulement s'il existe un isomorphisme $g \in End( \kq(-\gamma))$ tel que $f=f^{'}\circ g$}.

\begin{proof}
C'est la d\'efinition du sch\'ema $\QQuot_{\SS^{2(n-1)}\ku/\P^{n}}$.
 
\end{proof}

Le th\'eor\`eme suivant est une g\'en\'eralisation du th\'eor\`eme 4.1\cite{ca}.

\subsection{\bf Th\'eor\`eme} \label{4.8} 

{\em Soient $\kq$ le fibr\'e de quotient pond\'er\'e par les poids $\gamma,\alpha,\beta$ sur $\P^{n}$ et $\kf_{0}$ le~ fibr\'e de Tango pond\'er\'e par les poids $\gamma,\alpha,\beta$ sur $\P^{n}$. Les deux fibr\'es $\kq$ et $\kf_{0}$ sont d\'efinis par les suites exactes 
$$0\longrightarrow \kq(-\gamma) \stackrel{f_{0}}{\longrightarrow}\SS^{2(n-1)}\ku \stackrel{}{\longrightarrow}  \kf_{0}(\gamma) \longrightarrow 0$$

et
$$0\longrightarrow \ko_{\P^{n}}(-\gamma)\stackrel{}{\longrightarrow} \SS^{n}\ku\longrightarrow \kq \longrightarrow 0,$$

o\`u $f_{0}\in Hom(\kq(-\gamma),\SS^{2(n-1)}\ku)$. Alors chaque d\'eformation miniversale du fibr\'e $\kf_{0}$ est encore un fibr\'e de Tango pond\'er\'e sur $\P^{n}$. L'espace de Kuranishi de $\kf_{0}$ est lisse au point correspondant de $\kf_{0}$.}

\begin{proof}
Soient $f_{0}\in Hom(\kq(-\gamma),\SS^{2(n-1)}\ku) $ et $\kf_{0}(\gamma)=coker(f_{0})$ un fibr\'e \mbox{vectoriel} quotient de $\SS^{2(n-1)}\ku$ correspondant au morphisme $f_{0}$. Soit $Y \subseteq \QQuot_{\SS^{2(n-1)}\ku/\P^{n}}$ un \mbox{composant} irr\'eductible de $\QQuot_{\SS^{2(n-1)}\ku/\P^{n}}$ tel que  $f_{0} \in Y$. Soient $x \in Kur(\kq)$ \mbox{correspondant} au fibr\'e $\kq$ et \mbox{$z_{0} \in Kur(\kf_{0})$} correspondant au fibr\'e $\kf_{0}$.  
\xmat{  (Y,f_{0}) \ar[rr]^{\Psi}\ar[dd]_{\Phi} &&(Kur\kq,x)\\
\\ 
(Kur\kf_{0},z_{0})&&\\}

D'apr\`es le th\'eor\`eme \ref{4.5}, pour le morphisme $\Psi$, on a
$$  dim_{f_{0}}(Y)=dim_{x}(Kur(\kq))+ dim_{f_{0}}(\Psi^{-1}(x))$$

et $dim_{x}(Kur(\kq))= h^{1}(\EEnd(\kq))$. La dimension de la fibre du morphisme $\Psi$ est \'egale \`a $h^{0}(\kq^{*}(\gamma)\otimes \SS^{2(n-1)}\ku )-h^{0}(\EEnd(\kq))$, donc on obtient 
$$  dim_{f_{0}}(Y)=h^{1}(\EEnd(\kq))+ h^{0}(\kq^{*}(\gamma)\otimes \SS^{2(n-1)}\ku )-h^{0}(\EEnd(\kq)).$$

Pour le morphisme $\Phi$, d'apr\`es le th\'eor\`eme 3.12 page 137 et le th\'eor\`eme 2.2 page 126 \cite{qin}, on obtient 
$$h^{1}(\EEnd(\kf_{0}))\geq dim_{z_{0}}(Kur(\kf_{0}))\geq  dim_{f_{0}}(Y)- dim_{f_{0}}(\Phi^{-1}(z_{0})).$$

Soit 
$$Z=\lbrace f_{1}\in Y | \hspace{0.2 cm}\kf_{1}\simeq \kf_{0} \hspace{0.2cm} \text{o\`u $\kf_{1}$ est le fibr\'e correspondant \`a $f_{1}$ } \rbrace.$$

On obtient que
$$(\Phi^{-1}(z_{0}),f_{0})\subseteq (Z,f_{0})\hspace{0.2cm} et\hspace{0.2cm} dim_{f_{0}}((\Phi^{-1}(z_{0}),f_{0}))\leq dim_{f_{0}}((Z,f_{0})).$$

On en d\'eduit 
$$dim_{z_{0}}(Kur(\kf_{0}))\geq  dim_{f_{0}}(Y)-dim_{f_{0}}((Z,f_{0})).$$

Soit $\Sigma =\lbrace \sigma \in End(\SS^{2(n-1)}\ku)| \hspace{0.2 cm}\sigma .f_{0}=f_{0}\rbrace$. D'apr\`es les lemmes \ref{4.6} et \ref{4.7}, il en d\'ecoule que 
$$dim_{f_{0}}(Z)= h^{0}(\EEnd(\SS^{2(n-1)}\ku))-dim_{f_{0}}(\Sigma)-h^{0}(\EEnd(\kq)).$$

En consid\'erant la suite exacte suivante de fibr\'es vectoriels
$$0\longrightarrow \kf_{0}^{*}(-\gamma)\otimes \SS^{2(n-1)}\ku \longrightarrow \EEnd(\SS^{2(n-1)}\ku)\stackrel{}{\longrightarrow} \kq^{*}(\gamma)\otimes \SS^{2(n-1)}\ku  \longrightarrow 0,$$

on obtient la suite exacte suivante de groupes cohomologiques 
$$0\longrightarrow H^{0}(\kf_{0}^{*}(-\gamma)\otimes \SS^{2(n-1)}\ku) \longrightarrow End(\SS^{2(n-1)}\ku)\stackrel{\bullet\circ f_{0}}{\longrightarrow} H^{0}(\kq^{*}(\gamma)\otimes \SS^{2(n-1)}\ku) $$

qui nous donne $dim_{f_{0}}(\Sigma)=h^{0}(\kf_{0}^{*}(-\gamma)\otimes \SS^{2(n-1)}\ku)$. Donc on obtient que 
$$dim_{z_{0}}(Kur(\kf_{0}))\geq h^{1}(\EEnd(\kq))+h^{0}(\kq^{*}(\gamma)\otimes \SS^{2(n-1)}\ku)$$
$$- h^{0}(\EEnd(\SS^{2(n-1)}\ku))+h^{0}(\kf_{0}^{*}(-\gamma)\otimes \SS^{2(n-1)}\ku).$$

De la suite exacte pr\'ec\'edente de fibr\'es vectoriels, il s'ensuit que 
$$dim_{z_{0}}(Kur(\kf_{0}))\geq h^{1}(\EEnd(\kq))+ h^{1}(\kf_{0}^{*}(-\gamma)\otimes \SS^{2(n-1)}\ku). $$

En consid\'erant la suite exacte suivante de fibr\'es vectoriels
$$0\longrightarrow  \kf_{0}^{*}(-2\gamma)\otimes \kq  \longrightarrow  \kq(-\gamma)\otimes (\SS^{2(n-1)}\ku )^{*} \longrightarrow \EEnd(\kq) \longrightarrow 0$$

et comme on a $H^{1}(\kq(-\gamma)\otimes (\SS^{2(n-1)}\ku )^{*} )=H^{2}(\kq(-\gamma)\otimes (\SS^{2(n-1)}\ku )^{*} ) =0$, alors  
$$H^{1}(\EEnd(\kq))=H^{2}(\kf_{0}^{*}(-2\gamma)\otimes \kq) .$$

En consid\'erant la suite exacte suivante de fibr\'es vectoriels
$$0\longrightarrow  \kf_{0}^{*}(-2\gamma)\otimes \kq  \longrightarrow  \kf_{0}^{*}(-\gamma)\otimes \SS^{2(n-1)}\ku  \longrightarrow \EEnd(\kf_{0}) \longrightarrow 0,$$

on obtient la suite exacte suivante de groupes cohomologiques
$$\ldots \longrightarrow  H^{1}(\kf_{0}^{*}(-\gamma)\otimes \SS^{2(n-1)}\ku) \longrightarrow H^{1}(\EEnd(\kf_{0})) \longrightarrow H^{1}(\EEnd(\kq))\longrightarrow \ldots $$

qui nous donne $h^{1}(\EEnd(\kf_{0}))\leq h^{1}(\EEnd(\kq))+ h^{1}(\kf_{0}^{*}(-\gamma)\otimes \SS^{2(n-1)}\ku) $. Alors il en r\'esulte que \\ $dim_{z_{0}}(Kur(\kf_{0}))= h^{1}(\EEnd(\kf_{0}))$ et que $Kur(\kf_{0})$ est lisse en $z_{0}$. De plus on obtient
$$dim_{z_{0}}(Kur(\kf_{0}))= dim_{f_{0}}(Y)- dim_{f_{0}}(\Phi^{-1}(z_{0})).$$

D'apr\`es le th\'eor\`eme de la semi-continuit\'e des fibres (12.8, page 288 \cite{ha}), on en d\'eduit que\\ $dim_{z_{0}}(Im(\Phi))= dim_{z_{0}}(Kur(\kf_{0}))$ et que $\Phi $ est surjectif. Cela implique que $\kf_{0}$ est invariant par rapport \`a une d\'eformation miniversale.
 
\end{proof}

\end{document}